\documentclass[10pt,twoside]{article}
\usepackage[latin1]{inputenc}
\usepackage{amsmath, bm}
\usepackage{graphicx}
\usepackage{yhmath}
\usepackage[table]{xcolor}
\usepackage{mathrsfs} 
\usepackage{amssymb}
\usepackage{url}
\usepackage{makecell}
\usepackage{arydshln}
\usepackage{multirow}
\usepackage{arydshln}
\usepackage{mathtools}
\usepackage{stmaryrd}  

\input epsf
\def\limiten{\renewcommand{\arraystretch}{0.5}
\begin{array}[t]{c}\stackrel{}{\longrightarrow} \\
{\scriptstyle n\rightarrow
\infty}\end{array}\renewcommand{\arraystretch}{1}}

\def\limitepsn{\renewcommand{\arraystretch}{0.5}
\begin{array}[t]{c}\stackrel{a.s.}{\longrightarrow} \\
{\scriptstyle n \rightarrow
\infty}\end{array}\renewcommand{\arraystretch}{1}}

\def\limiteloin{\renewcommand{\arraystretch}{0.5}
\begin{array}[t]{c}\stackrel{{\cal D}}{\longrightarrow} \\
{\scriptstyle n\rightarrow
\infty}\end{array}\renewcommand{\arraystretch}{1}}

\numberwithin{equation}{section}

\newtheorem{thm}{Theorem}[section]

\newtheorem{lem}[thm]{Lemma}

\newtheorem{prop}[thm]{Proposition}

\newcommand{\E}{\ensuremath{\mathbb{E}}}
\newcommand{\R}{\ensuremath{\mathbb{R}}}
\newcommand{\Z}{\ensuremath{\mathbb{Z}}}

\newcommand{\N}{\ensuremath{\mathbb{N}}}

\definecolor{grisclair}{gray}{0.9}
\font\dsrom=dsrom10 scaled 1200

\renewcommand{\arraystretch}{.8}

\begin{document}
\title{\bf Density power divergence for general integer-valued time series with  exogenous covariates }
 \maketitle \vspace{-1.0cm}
 \begin{center}
   Mamadou Lamine DIOP \footnote{Supported by
   the MME-DII center of excellence (ANR-11-LABEX-0023-01) 
   } 
   and 
    William KENGNE \footnote{Developed within the ANR BREAKRISK: ANR-17-CE26-0001-01 and the  CY Initiative of Excellence (grant "Investissements d'Avenir" ANR-16-IDEX-0008), Project "EcoDep" PSI-AAP2020-0000000013} 
 \end{center}

  \begin{center}
  { \it 
 THEMA, CY Cergy Paris Université, 33 Boulevard du Port, 95011 Cergy-Pontoise Cedex, France.\\
  E-mail: mamadou-lamine.diop@u-cergy.fr ; william.kengne@u-cergy.fr  \\
  }
\end{center}

 \pagestyle{myheadings}
 \markboth{MDPDE for integer-valued time series}{Diop and Kengne}

~~\\
\textbf{Abstract}:
 In this article, we study a robust estimation method for a general class of integer-valued time series models. 
 The conditional distribution of the process belongs to a broad class of distribution and unlike classical autoregressive framework, the conditional mean of the process also depends on some exogenous covariates.      
 We derive a robust inference procedure based on the minimum density power divergence. 
  Under certain regularity conditions, we establish that the proposed estimator is consistent and asymptotically normal. 
  In the case where the conditional distribution belongs to the exponential family, we provide sufficient conditions for the existence of a stationary and ergodic $\tau$-weakly dependent solution.
 Simulation experiments are conducted to illustrate the empirical performances of the estimator. An application to the number of transactions per minute for the stock
Ericsson B is also provided.\\ 


\section{Introduction}
The analysis of time series of counts  has attracted much interest
in the literature during the last two decades, 
 given the large number of papers written in this direction. 
  This is due among others to the applications in various fields: 
  epidemiological surveillance (number of new infections), in finance (number of transactions), in industrial quality control (number of defects), and traffic accidents (number of road casualties), etc. 
 To describe time series of count data, several questions have been addressed with various modeling approaches which are generally classified into two categories: observation-driven models and parameter-driven models (see Cox (1981)).
 One of the first important results on this topic were obtained independently by McKenzie (1985) and  Al-Osh and Alzaid (1987); the INAR model has been introduced by using the binomial thinning operator. Due to its limitations, numerous extensions have been proposed; see, e.g., Al-Osh and Alzaid (1990).  
 Later, new models with various marginal distributions and dependence structures have been studied by several authors; see among others, Fokianos {\it et al.} (2009), Doukhan \textit{et al.} (2012, 2013), Doukhan and Kengne (2015), Davis and Liu (2016), Ahmad and Francq (2016), Douc \textit{et al.} (2017),  Fokianos {\it et al.} (2020) for some recent progress.   
%

\medskip 

 For most developments in the literature, the parametric inference is commonly based on the conditional maximum likelihood estimator (MLE) which provides a full asymptotic efficiency among regular estimators.
  However, it has been recognized that the MLE is very sensitive to small perturbations caused by outliers in the underlying model.
Beran (1977) addressed this issue and was one of the first references in the literature to use the density-based minimum divergence methods.  Numerous others works have been devoted to this topic; see, among others, Tamura and Boos (1986), Simpson (1987), Basu and Lindsay (1994), and Basu \textit{et al.} (1998). 
 In the context of modelling time series of counts, this question has already been investigated. 
%
For instance, Fokianos and Fried (2010, 2012) have studied the problem of intervention effects (that generating various types of outliers) in linear and  log-linear Poisson autoregressive models.
Fried {\it et al.} (2015) have proposed a Bayesian approach for handling additive outliers in INGARCH processes; they have used Metropolis-Hastings algorithm to estimate the parameters of the model.
 Recently, Kim and Lee (2017, 2019)  adopted the approach of Basu \textit{et al.} (1998) to construct a robust estimator based on the  minimum density power divergence for zero-inflated Poisson
autoregressive models, and a general integer-valued time series whose conditional distribution belongs to the one-parameter exponential family.
See also Kang and Lee (2014) for application of such procedure to Poisson autoregressive models.
   
\medskip

On the other hand, most of the models proposed for analysing  time series of counts do not provide a framework within which we are able to analyze the possible dependence of the observations on a relevant exogenous covariates. 
%
 In this vein, Davis and Wu (2009) have studied generalized linear models for time series of counts, where conditional on covariates, the observed process is modelled by a negative binomial distribution.
 Agosto {\it et al.} (2016) have developed a class of linear Poisson autoregressive models with exogenous covariates (PARX), 
where the parametric inference is based on the maximum likelihood method.
 Later, Pedersen and Rahbek (2018) proposed a theory for testing the significance of covariates in a class of PARX models.
 See also the recent work of Fokianos and Truquet (2019) which considered   a class of categorical time series models with covariates and addressed stationarity and ergodicity question.
 The R package "tscount"  Liboschik {\it et al.} (2017) provides tools for  analysis and modeling count time series models where the conditional mean
of the process can take into account covariate effects.

\medskip
 
In this new contribution, we consider a quite general class of observation-driven models for time series of counts with  exogenous covariates. 
 The mean of the discrete conditional distribution of $Y_t$ depends on the whole past observations, some relevant covariates and a finite dimensional parameter $\theta^*$. 
 We study a robust estimator of $\theta^*$ by using the minimum density power divergence estimator (MDPDE) proposed by Basu \textit{et al.} (1998). 
 Compared to that of Kim and Lee (2019), the framework considered here is more general: 
   (i) the class of models has the ability to handle the dependence of the observations on exogenous covariates, 
  (ii) the dependence through the past is of infinite order (which enables a large dependence structure and to consider INGARCH($p,q$)-type models) and 
  (iii) even if in many applications the conditional distribution belongs to the exponential family, the class of the conditional distribution considered here is beyond the exponential family. 
  \medskip
 
 Recently, Aknouche A. and Francq (2020) consider a class of count time series 
 with the conditional distributions that are not restricted to belong to the one-parameter exponential family and where the conditional mean depend on
exogenous covariates. They provided conditions for  stationarity and ergodicity.
  For both the linear and nonlinear conditional means considered by these authors, the feedback effects of the covariate is linear and it is an additive term in the conditional mean.
  We do not restrict to such setting and consider a larger dependence structure between the conditional means and the covariate. In the case where the conditional distribution belongs to the exponential family, we provide sufficient conditions that ensure the existence of a stationary and ergodic $\tau$-weakly dependent solution. In this sense of a large dependence to the covariate, our result is more general.   
  %
  %

\medskip

The paper is structured as follows. Section 2 contains the model specification and the construction of the robust estimator 
as well as the main results. Section 3 is devoted to the application of the general results to some examples of dynamic models. 
Some simulation results are displayed in Section 4 whereas Section 5 focus on applications on a real data example.
The proofs of the main results are provided in Section 6.



 \section{Model specification and estimation}
 %
 %
 \subsection{Model formulation} 
  Suppose that $\{Y_{t},\,t\in \Z \}$ is a time series of counts and that $X_t=(X_{1,t},X_{2,t},\cdots,X_{d_{x}, t}) \in \R^{d_{x}}$ represents a vector of covariates with $d_{x} \in \N$. 
  Denote by $\mathcal{F}_{t-1}=\sigma\left\{Y_{t-1},\cdots ;X_{t-1},\cdots  \right\}$ the $\sigma$-field generated by the whole past at time $t-1$. 
  Consider the general dynamic model where $Y_t | \mathcal{F}_{t-1}$ follows a discrete distribution whose mean satisfying:
  
   \begin{equation}\label{Model}
     \E(Y_t | \mathcal{F}_{t-1}) =\lambda_t(\theta) 
      = f_\theta(Y_{t-1}, Y_{t-2}, \cdots; X_{t-1},X_{t-2}\cdots)
   \end{equation}
  where 
   $f_\theta$ is a measurable non-negative function defined on $\mathbb{N}_0^{\N} \times (\R^{d_{x}})^{\N}$ (with $\N_0 = \N\cup \{ 0\}$) and assumed to be known up to the parameter $\theta$ which belongs in a compact subset $\Theta \subset \R^d$ ($d \in \N$).
    Let us note that when $X_t \equiv C$ (a constant), then the model (\ref{Model}) reduces to the classical integer-valued autoregressive model that has already been considered in the literature (see, e.g., Ahmad and Francq (2016)). 
    In the sequel, we assume that the random variables $Y_t,~ t \in \Z$ have the same distribution and denote  by $G_\theta(\cdot|\mathcal{F}_{t-1})$ the distribution of $Y_t|\mathcal{F}_{t-1}$; let $g (\cdot|\eta_t)$ be the probability density function of this distribution, where  $\eta_t$ is  the natural parameter of $G_\theta$ given by $\eta_t=\eta(\lambda_t(\theta))$.
     Assume that the density $g (\cdot|\eta_t(\theta))$ is known up to the parameter $\theta$; and this density has a support set $\{y,~ g (y|\eta_t(\theta))>0\}$ which is independent of $\theta$. 
  
  \medskip
  
\noindent
 Throughout the sequel, the following norms will be used:
 {\em
\begin{itemize}
\item $\left\|x\right\| \coloneqq  \sup_{1  \leq i  \leq p} |x_i|$, for any $x \in \mathbb{R}^{p}$, $p \in \N$;
\item  $\left\|f\right\|_{\Theta} \coloneqq \sup_{\theta \in \Theta}\left(\left\|f(\theta)\right\|\right)$ for any function $f:\Theta \longrightarrow   M_{p,q}(\R)$,
 where $M_{p,q}(\R)$ denotes the set of matrices of dimension $p\times q$ with coefficients in $\R$, for $p,q \in \N$;
\item $\left\|Y\right\|_r \coloneqq \E\left(\left\|Y\right\|^r\right)^{1/r}$, if $Y$ is a random vector with finite $r-$order moments, for $r >0$.
\end{itemize}
}

 \medskip
 
\noindent
We set the following classical Lipschitz-type condition on the function $f_\theta$.

    \medskip
    \noindent \textbf{Assumption} \textbf{A}$_i (\Theta)$ ($i=0,1,2$):
    For any $(y,x) \in \mathbb{N}_0^{\N} \times \R^{\infty}$, the function $\theta \mapsto f_\theta(y,x)$ is $i$ times continuously differentiable on $\Theta$  with $ \left\| \partial^i f_\theta(0)/ \partial \theta^i\right\|_\Theta<\infty $; 
    and
      there exists two sequences of non-negative real numbers $(\alpha^{(i)}_{k,Y})_{k\geq 1} $ and $(\alpha^{(i)}_{k,X})_{k\geq 1} $ satisfying
     $ \sum\limits_{k=1}^{\infty} \alpha^{(0)}_{k,Y} <1 $ (or $ \sum\limits_{k=1}^{\infty} \alpha^{(i)}_{k,Y} <\infty $ for $i=1, 2$) and $ \sum\limits_{k=1}^{\infty} \alpha^{(i)}_{k,X} <\infty $ for $i=0,1, 2$;
   such that for any  $(y,x), (y',x') \in \mathbb{N}_0^{\N} \times (\R^{d_{x}})^{\N}$,
  \[  \Big \| \frac{\partial^i f_\theta(y,x)}{ \partial \theta^i}-\frac{\partial^i f_\theta(y',x')}{\partial\theta^i} \Big \|_{\Theta}
  \leq  \sum\limits_{k=1}^{\infty}\alpha^{(i)}_{k,Y}  |y_k-y'_k| +  \sum\limits_{k=1}^{\infty}\alpha^{(i)}_{k,X} \left\|x_k-x'_k \right\|. 
  \]
where $\| \cdot\|$ denotes any vector or matrix norm.

 \medskip 
 
 \noindent
 In the sequel, it is assumed for the general setting that there exists a stationary and ergodic process $Y^*_{t}=(Y_{t},\lambda_t, X_{t})$ solution of (\ref{Model}); and 
 \begin{equation}\label{moment}
    \exists C, \epsilon >0, \text{ such that } \forall t \in \Z, ~ ~  \E \left\|Y^*_{t}\right\|^{1+\epsilon} < C.
   \end{equation}
 In the case where the distribution of $Y_t | \mathcal{F}_{t-1}$ belongs to the one-parameter exponential family, the existence of such solution is established (see Proposition \ref{exp_existence}). 
  
 \subsection{Minimum power divergence estimator}
 In this subsection, we briefly describe the use of the density power divergence to obtain an estimation of the parameters of the model (\ref{Model}). The asymptotic behavior of the estimated parameter is also studied.    
Assume that the observations $(Y_1,X_1),\cdots,(Y_n,X_n)$ are generated from (\ref{Model}) according to the true parameter $\theta^* \in \Theta$ which is unknown; i.e., $g (\cdot|\eta_t(\theta^*))$ is the true conditional density of $Y_t|\mathcal{F}_{t-1}$.  
 Let $\mathbb G=\left\{ g(\cdot|\eta_t(\theta)); ~ \theta \in \Theta \right\}$ be the parametric family of density functions indexed by $\theta \in \Theta$.
 To estimate $\theta^*$, Basu \textit{et al.} (1998) have proposed a method which consists to choose
  the "best approximating distribution" of $Y_t|\mathcal{F}_{t-1}$ in the family $\mathbb G$ by minimizing  
  a divergence $d_\alpha$ between density functions $g(\cdot|\eta_t(\theta))$ and $g(\cdot|\eta_t(\theta^*))$.
The density power divergence $d_\alpha$ between two density functions $g$ and $g_{*}$ is defined by (in the discrete set-up)
  \[
 d_\alpha(g,g_{*})=\left\{
\begin{array}{ll}
   \sum\limits_{y = 0}^{\infty} \left\{g^{1+\alpha}(y)  -\big( 1+\frac{1}{\alpha}\big)
   g_{*}(y) g^{\alpha}(y) + \frac{1}{\alpha}g^{1+\alpha}_{*}(y)\right\}, & \alpha >0 , 
    \\
   \\
   \sum\limits_{y = 0}^{\infty}\big\{ g_{*}(y) \big(\log g_{*}(y) -\log g(y) \big)\big\} , & \alpha =0. \\
\end{array}
\right.
 \]
 So, the empirical objective function (based on the divergence between conditional density functions) up to some terms which are independent of $\theta$ is
 $H_{\alpha,n}(\theta)=\frac{1}{n}  \sum_{t =1}^{n} \ell_{\alpha,t}(\theta)$ where 
 \[
 \ell_{\alpha,t}(\theta)=\left\{
\begin{array}{ll}
   \sum\limits_{y = 0}^{\infty}g (y|\eta_t(\theta))^{1+\alpha}   -\big( 1+\frac{1}{\alpha}\big)g (Y_t|\eta_t(\theta))^{\alpha},& \alpha >0 ,  \\
  
   \\
  - \log g (Y_t|\eta_t(\theta)) ,& \alpha =0, \\
\end{array}
\right.
 \]
with $\eta_t(\theta)= \eta(\lambda_t(\theta) )$ and $\lambda_t(\theta) 
      = f_\theta(Y_{t-1}, Y_{t-2}, \cdots; X_{t-1},X_{t-2}\cdots)$. 
 %
 %
Since $(Y_0,X_0),(Y_{-1},X_{-1}),\cdots$ are not observed, $H_{\alpha,n}(\theta)$ is approximated by 
$\widehat H_{\alpha,n}(\theta) =\frac{1}{n}\sum\limits_{t =1}^{n} \widehat \ell_{\alpha,t}(\theta)$, where
 
 \[
\widehat  \ell_{\alpha,t}(\theta)=\left\{
\begin{array}{ll}
   \sum\limits_{y = 0}^{\infty} g (y|\widehat \eta_t(\theta))^{1+\alpha}   -\big( 1+\frac{1}{\alpha}\big)g (Y_t|\widehat \eta_t(\theta))^{\alpha},& \alpha >0 ,  \\
  
   \\
  - \log g (Y_t|\widehat \eta_t(\theta)) ,& \alpha =0, \\
\end{array}
\right.
 \]
with $\widehat \eta_t(\theta)= \eta(\widehat \lambda_t(\theta) )$ and $\widehat \lambda_t(\theta) 
      = f_\theta(Y_{t-1}, \cdots,Y_1,0,\cdots; X_{t-1},\cdots,X_1,0,\cdots)$.
 Therefore, the MDPDE of $\theta^*$ is defined by (cf. Basu \textit{et al.} (1998) and Kim and Lee (2019))
   \[
 \widehat{\theta}_{\alpha,n}= \underset{\theta\in \Theta}{\text{argmin}} ( \widehat H_{\alpha,n}(\theta) ).
 \]
 Let us recall that when $\alpha=0$, the MDPDE corresponds to the MLE. 

        \medskip 
          
 \noindent
  We need the following regularity assumptions to study the consistency and the asymptotic normality of the MDPDE. 
 \begin{enumerate}
   \item [(\textbf{A0}):]  for all  $\theta\in \Theta$,
 $ \big( f_\theta(Y_{t-1},Y_{t-2}, \cdots ; X_{t-1},X_{t-2},\cdots)= f_{\theta^*}(Y_{t-1},Y_{t-2}, \cdots ; X_{t-1},X_{t-2},\cdots)  \ \text{a.s.} ~ \text{ for some } t \in \Z \big) \Rightarrow ~ \theta = \theta^*$;
  moreover, $\exists  \underline{c}>0$ such that $\displaystyle \inf_{ \theta \in \Theta} f_\theta(y_{1}, \cdots ; x_{1},\cdots)  \geq \underline{c}$ for all $ (y,x) \in  \N_0^{\N} \times \R^\infty$ ;
  \item [(\textbf{A1}):] $\theta^*$ is an interior point in the compact parameter space $\Theta \subset \mathbb{R}^{d}$;
  \item [(\textbf{A2}):]  There exists a constant $\underline{\lambda}>0$ such that $~\E \left\|\lambda_t(\theta)\right\|^4_\Theta =\underline{\lambda}<\infty$, for all $t \geq 1$;
  \item [(\textbf{A3}):] for all $y \in \N_0$, the function $\eta \mapsto g(y|\eta)$ is twice continuously differentiable on $\R$ and for some $\eta, \eta' \in \R$,
  $\big( g(y|\eta) = g(y|\eta') ~ \forall y \in \N_0   \big) \Rightarrow \eta = \eta'$;
  \item [(\textbf{A4}):] the mapping $\eta \mapsto\varphi(\eta)= \sum_{y =0}^{\infty} \left| \frac{\partial g (y| \eta)}{\partial \eta}\right|$ is well definite on $\R$ and for all $t \geq 1$, there exists a constant $\overline{\varphi}_t>0 $ such that 
         $  \underset{0 \leq \delta \leq 1}{\sup} \left \|  \left \| \varphi \left ( \delta \eta_t(\theta) + (1-\delta) \widehat{\eta}_t(\theta)\right) \right\|_{\Theta}\right \|_2 \leq \overline{\varphi}_t < \infty$;
   \item [(\textbf{A5}):] for all $t \geq 1$, the mapping $\eta \mapsto \psi(\eta)= \left|\frac{1}{g (Y_t|  \eta)} \frac{\partial g (Y_t|  \eta)}{\partial \eta}\right|$ is well definite on $\R$ and  there exists a constant $\overline{\psi}_t>0 $ such that 
         $  \underset{0 \leq \delta \leq 1}{\sup} \left \|  \left \| \psi \left ( \delta \eta_t(\theta) + (1-\delta) \widehat{\eta}_t(\theta)\right) \right\|_{\Theta}\right \|_2  \leq \overline{\psi}_t < \infty$;  
  \item [(\textbf{A6}):] the mapping $\lambda \mapsto \eta(\lambda)$ (defined in $[\underline{c},\, +\infty[$) satisfying the Lipschitz condition: there exists a constant $c_{\eta} >0$ such that, for all $\lambda, \lambda' >0$, $|\eta(\lambda) - \eta(\lambda')| \leq c_{\eta} |\lambda - \lambda'|$; moreover, $\big(\eta(\lambda) = \eta(\lambda')  \big) \Rightarrow \lambda = \lambda'$.
    \item [(\textbf{A7}):]  there exists two non-negative constants $\overline{h}_{\alpha}$ and $\overline{m}_{\alpha}$  such that the mappings: 
    \begin{align}
    \label{def_h_alpha}&\lambda \mapsto h_\alpha(\lambda)=(1+\alpha)\eta'(\lambda)\bigg[\underset{y =0}{\overset{\infty}{\sum}} \frac{\partial g (y|\eta(\lambda))}{ \partial \eta}g (y|\eta(\lambda))^{\alpha}  -\frac{\partial g (Y_t|\eta(\lambda))}{ \partial \eta}g (Y_t|\eta(\lambda))^{\alpha-1} \bigg] \\
     \text{and}~~~~ &\nonumber\\ 
     \label{def_m_alpha}&\lambda \mapsto m_\alpha(\lambda)=\frac{\partial h_\alpha(\lambda)}{\partial \lambda}
    \end{align}
      satisfy    
    $\big\|\left\|h_\alpha(\lambda_t(\theta)) \right\|^2_\Theta\big\|_2 \leq \overline{h}_{\alpha}~$ 
     and 
    $~\underset{0 \leq \delta \leq 1}{\sup}  \big\|\| m_\alpha \left( \delta \lambda_t(\theta) + (1-\delta) \widehat\lambda_t(\theta)\right) \|_{\Theta}\big\|_2  \leq \overline{m}_{\alpha}$, for all $t \geq 1$;
    \item [(\textbf{A8}):] for all $c^T \in \R$, $c^T \frac{\partial \lambda_{t} (\theta^*)}{\partial    \theta}=0$ a.s   $\Longrightarrow~ c^T=0$,  where $^T$ denotes the transpose.
\end{enumerate}

  \medskip  
   
 \noindent Assumption (\textbf{A0}) is an identifiability condition.  The conditions (\textbf{A3})-(\textbf{A7}) allow us to unify the theory and the treatment for a class of distributions that belongs or not to the exponential family. In the class of the exponential family distribution, these conditions can be written in a simpler form (see the Subsection \ref{Sec_Exempl_Ex}). The other conditions are standard in this framework and can also be found in many studies; see for instance,  Kim and Lee (2019). As detailed in Section \ref{Sec_Example}, all these conditions are satisfied for many classical models.
 
  \medskip  

 \noindent
 The following theorem gives the consistency of the estimator $\widehat{\theta}_{\alpha,n}$.

\begin{thm}\label{th1}
Assume that (\textbf{A0})-(\textbf{A6}),  (\textbf{A}$_0(\Theta)$) and (\ref{moment}) (with $\epsilon > 2$) hold with 
\begin{equation}\label{eq_th1}
 \alpha^{(0)}_{k,Y} +\alpha^{(0)}_{k,X}= \mathcal{O}(k^{-\gamma}) \text{ for some } \gamma >1; \text{ and } 
 \sum_{k \geq 1} \frac{1}{k^{\gamma}} \max \left( \overline{\varphi}_k, \overline{\psi}_k\right) < \infty.
\end{equation}
Then  
\begin{equation*}\label{Cons_theta}
\widehat{\theta}_{\alpha,n} \limitepsn \theta^*.
\end{equation*}
\end{thm}
~\\
 The following theorem gives the asymptotic normality of $\widehat{\theta}_{\alpha,n}$.

\begin{thm}\label{th2}
Assume that (\textbf{A0})-(\textbf{A8}),  (\textbf{A}$_i(\Theta)$) (for $i=0,1,2$), (\ref{moment}) (with $\epsilon > 3$) 
 with

\begin{equation}\label{eq_th2}
  \alpha^{(0)}_{k,Y} +\alpha^{(0)}_{k,X} +\alpha^{(1)}_{k,Y} +\alpha^{(1)}_{k,X} = \mathcal{O}(k^{-\gamma}) \text{ for some } \gamma >3/2.
\end{equation}
 Then  

\begin{equation*}\label{Normal_theta}
 \sqrt{n}\left(\widehat\theta_{\alpha,n}-\theta^*\right) \limiteloin \mathcal{N}_d \left(0,\Sigma_\alpha \right) \ \text{ with } \
 \Sigma_\alpha=J^{-1}_\alpha I_\alpha J^{-1}_\alpha, 
\end{equation*}
 where 
 $
 J_\alpha=- \E\Big(\frac{ \partial^2  \ell_{\alpha,1}(\theta^*)}{\partial \theta \partial \theta^T}\Big)
 \text{ and }
  I_\alpha= \E\Big[ \Big(\frac{ \partial  \ell_{\alpha,1}(\theta^*)}  {\partial \theta} \Big) \Big( \frac{ \partial  \ell_{\alpha,1}(\theta^*)}{\partial \theta} \Big)^T \Big].
 $
 \end{thm}
The trade-off between the robustness and the efficiency is controlled by the tuning parameter $\alpha$. As pointed out in numerous works (see for instance, Basu \textit{et al.} (1998)), it is found that the estimators with large $\alpha$ have strong robustness properties while small value of $\alpha$   is suitable when the efficiency is preferred. 
So, the procedure is typically less efficient when $\alpha$ increases.
In the empirical studies, we will consider the value of $\alpha$ between zero and one, and models with condition distribution belonging to the one-parameter exponential family.
%
 %


\section{Examples}\label{Sec_Example}
  In this section, we give some particular cases of class of integer-valued time series defined in (\ref{Model}). We show that the regularity conditions required for the asymptotic results of the previous section are satisfied for these models. 
  Throughout the sequel, we consider that $X_t=(X_{1,t},X_{2,t},\cdots,X_{d_{x}, t}) \in \R^{d_{x}}$ ($d_{x} \in \N$) represents  a vector of covariates, 
  $\theta$ belongs to a compact set $\Theta \subset  \mathbb{R}^{d}$ ($d \in \N$) and $C$ denotes a positive constant whom value may differ from one inequality to another. 
  For any $\theta \in \Theta$, we will use the notation  $\eta_{t,\delta}(\theta) \coloneqq  \delta \eta_t(\theta) + (1-\delta) \widehat{\eta}_t(\theta)$ (with $0\leq \delta \leq 1$) for a given element between $\widehat \eta_{t}(\theta)$ and $\eta_{t}(\theta)$.
  
  \subsection{A general model with the exponential family distribution}\label{Sec_Exempl_Ex}

   As a first example, we consider a process $\{Y_{t},\, t\in \mathbb{Z}\}$ satisfying:
\begin{equation} \label{Model_Expo}
Y_t|\mathcal{F}_{t-1}\sim g(y|\eta_t)~~\textrm{with}~~\lambda_t \coloneqq \E(Y_t|\mathcal{F}_{t-1})=f_{\theta}(Y_{t-1}, \cdots; X_{t-1},\cdots)
\end{equation}
where 
  $g(\cdot|\cdot)$ is a discrete distribution that belongs to the one-parameter exponential family; that is,
\[ g(y|\eta)=\exp\left\{\eta y -A(\eta)\right\}h(y) \]
where $\eta$ is the natural parameter (i.e., $\eta_t$ is the natural parameter of the distribution of $Y_t|\mathcal{F}_{t-1}$), $A(\eta)$, $h(y)$ are known functions and $f_{\theta}(\cdot)$ is a non-negative function defined on
$\N^{\infty}_0 \times \R^{\infty}$, assumed to be know up to $\theta$. 
Let us set $B(\eta) =A^{\prime}(\eta)$ the  derivative of $A(\eta)$ (which is assumed to exist a well as the second order derivative); it is known that $\E(Y_t|\mathcal{F}_{t-1})=B(\eta_t)$. Therefore, the  
model (\ref{Model_Expo}) can be seen as a particular case of (\ref{Model}), where $\eta(\lambda)=B^{-1}(\lambda)$. 
Similar models  without covariates have been studied by Davis and Liu (2016) (with the MLE) and Kim and Lee (2019) (with the MDPDE) where the conditional mean of $Y_t$ depends  only on $(Y_{t-1},\lambda_{t-1})$.
Cui and Zheng (2017) carried out the model (\ref{Model_Expo}) without covariates and  proved the existence of a stationary, ergodic and $\tau$-weakly dependent solution; they also considered inference based on the  conditional maximum likelihood estimator.
 The following proposition establishes the existence of a unique solution for the class of model (\ref{Model_Expo}) with covariates. 
 Let us impose a autoregressive-type structure on the covariates:
 \begin{equation}\label{exp_Markov_cov}
  X_t=u(X_{t-1},X_{t-2},\cdots, \varepsilon_t),
 \end{equation}
 where ($\varepsilon_t$) is a sequence of independent and identically distributed (\textit{i.i.d}) random variables and $u(x; \varepsilon_t)$ a function with values in $\R^{d_x}$, satisfying 
  \begin{equation} \label{exp_Lip_cov}
  \E \| u(0; \varepsilon_t) \| < \infty  ~ \text{ and }  ~   \E \| u(x; \varepsilon_t) - u(x'; \varepsilon_t) \| \leq \sum_{k \geq 1} \alpha_k(u) \|x_k - x'_k \|  ~ \text{ for all } x, x'  \in (\R^{d_{x}})^{\N},
 \end{equation}  
 for some non-negative sequence $(\alpha_k(u))_{k \geq 1}$ such that  $\sum_{k \geq 1} \max\big\{ \alpha_k(u),  \alpha_{k,Y}^{(0)} \big\} < 1 $. 
 

 \begin{prop}\label{exp_existence}
 Assume that \textbf{A}$_0(\Theta)$ and (\ref{exp_Lip_cov}) holds. Then there exists a $\tau$-weakly dependent stationary solution $Y^*_t=(Y_t, \lambda_t, X_t)$ of (\ref{Model_Expo}) and (\ref{exp_Markov_cov}) such that $\E\| Y^*_t \| < \infty$.
 \end{prop}

\noindent 
For the model (\ref{Model_Expo}), let us provide some sufficient conditions for the assumptions (\textbf{A2})-(\textbf{A7}).
%

%
                                               
\begin{itemize}
\item  
According to Assumption \textbf{A}$_0(\Theta)$, for all $t \geq 1$, we have
\begin{align} 
 \|\lambda_{t}(\theta)\|_\Theta 
     &  \leq
     \left\| f_\theta(0)\right\|_\Theta
     +
     \left\| f_\theta(0) - f_\theta  ( Y_{t-1},\cdots ;X_{t-1},\cdots ) \right\|_\Theta \nonumber \\
    \label{App_A0} &  \leq
    \left\| f_\theta(0)\right\|_\Theta+   \sum\limits_{\ell \geq 1}  \alpha^{(0)}_{\ell,Y}  Y_{t-\ell}+\sum\limits_{\ell \geq 1}  \alpha^{(0)}_{\ell,X}\left\|  X_{t-\ell}\right\|.
    \end{align} 
    Thus, 
    \begin{align*} 
 \left\|\|\lambda_{t}(\theta)\|_\Theta\right\|_4 
    &  \leq
    \left\|\left\| f_\theta(0)\right\|_\Theta\right\|_4+   \sum\limits_{\ell \geq 1}  \alpha^{(0)}_{\ell,Y}  \left\|Y_{t-\ell}\right\|_4+ \sum\limits_{\ell \geq 1} \alpha^{(0)}_{\ell,X} \left\|\left\|  X_{t-\ell}\right\|\right\|_4\\
    &  \leq
   C+   C\big(\sum\limits_{\ell \geq 1}  \alpha^{(0)}_{\ell,Y} +\sum\limits_{\ell \geq 1}  \alpha^{(0)}_{\ell,X}\big)<\infty ~~(\text{from the assumption (\ref{moment}) with } \epsilon>3).
    \end{align*} 
    Hence, (\textbf{A2}) is satisfied.
  \item Clearly, (\textbf{A3}) is satisfied.
	\item  
	According to the above notations, for any $\theta \in \Theta$, $\delta \in [0,1]$, we have 
\begin{align*} 
\varphi(\eta_{t,\delta}(\theta))&= \sum_{y =0}^{\infty} \left| (y -B(\eta_{t,\delta}(\theta))) g(y| \eta_{t,\delta}(\theta))\right|
\leq 2 B(\eta_{t,\delta}(\theta)). 
\end{align*}
Since the function $B$ is strictly increasing (because  Var$(Y_t|\mathcal{F}_{t-1})=B'(\eta_t)>0$), we deduce that
\begin{equation}\label{phi_eta_delta_norm}
\left\|\varphi(\eta_{t,\delta}(\theta))\right\|_\Theta \leq
 2 \left( \left\|B(\eta_{t}(\theta))\right\|_\Theta+\left\|B(\widehat \eta_{t}(\theta))-B(\eta_{t}(\theta))\right\|_\Theta\right) =2 \left( \|\lambda_t(\theta)\|_\Theta+\|\widehat \lambda_{t}(\theta)-\lambda_{t}(\theta)\|_\Theta\right).
\end{equation}  

Moreover, Assumption \textbf{A}$_0(\Theta)$ implies 
\begin{align}
    \|\widehat \lambda_{t}(\theta)-\lambda_{t}(\theta)\|_\Theta 
    \nonumber &\leq  \left\| f_\theta ( Y_{t-1},\cdots,Y_{1},0,\cdots;X_{t-1},\cdots,X_{1},0,\cdots) - f_\theta  ( Y_{t-1},\cdots ;X_{t-1},\cdots ) \right\|_\Theta \\
     \label{App1_A0} &  \leq  \sum\limits_{\ell \geq t}  \alpha^{(0)}_{\ell,Y} Y_{t-\ell}+ \sum\limits_{\ell \geq t}  \alpha^{(0)}_{\ell,X}\left\|X_{t-\ell}\right\|\\
      \label{App1_prime_A0} &  \leq  \sum\limits_{\ell \geq 1}  \alpha^{(0)}_{\ell,Y} Y_{t-\ell}+ \sum\limits_{\ell \geq 1}  \alpha^{(0)}_{\ell,X}\left\|X_{t-\ell}\right\|.
      \end{align} 
      Thus, from (\ref{App_A0}), (\ref{phi_eta_delta_norm}) and (\ref{App1_prime_A0}), for $t \geq 1$, we get
      \begin{align*}
      \left \|  \left \| \varphi(\eta_{t,\delta}(\theta)) \right\|_{\Theta}\right \|_2 
      &\leq     
       2 \Big( \left\|\left\| f_\theta(0)\right\|_\Theta\right\|_2+  2\sum\limits_{\ell \geq 1}  \alpha^{(0)}_{\ell,Y} \|Y_{t-\ell}\|_2 +2\sum\limits_{\ell \geq 1}  \alpha^{(0)}_{\ell,X}\left\|  \|X_{t-\ell}\right\|\|_2 \Big)\\
       &\leq  
       C+  C\big(\sum\limits_{\ell \geq 1}  \alpha^{(0)}_{\ell,Y} +\sum\limits_{\ell \geq 1}  \alpha^{(0)}_{\ell,X}\big) <\infty ~(\text{from the stationary assumption}).
      \end{align*}  
 Therefore, (\textbf{A4}) is satisfied with $ \overline{\varphi}_t =  C+  C\big(\sum\limits_{\ell \geq 1}  \alpha^{(0)}_{\ell,Y} +\sum\limits_{\ell \geq 1}  \alpha^{(0)}_{\ell,X}\big)$ is constant. 
 In addition, remark that
 \begin{equation*}\label{psi_eta_delta_norm}
 \left\|\psi(\eta_{t,\delta}(\theta))\right\|_\Theta
 = \left\| (Y_t -B(\eta_{t,\delta}(\theta)))\right\|_\Theta \leq Y_t + \|\lambda_t(\theta)\|_\Theta+\|\widehat \lambda_{t}(\theta)-\lambda_{t}(\theta)\|_\Theta.
\end{equation*}  
Thus, from  (\ref{App_A0}), (\ref{App1_prime_A0}), we deduce 
\begin{align*} 
\left\|\left\|\psi(\eta_{t,\delta}(\theta))\right\|_\Theta\right\|_2
&\leq 
        \|Y_{t}\|_2+  \left\|\left\| f_\theta(0)\right\|_\Theta\right\|_2+  2\sum\limits_{\ell \geq 1}  \alpha^{(0)}_{\ell,Y} \|Y_{t-\ell}\|_2 +2\sum\limits_{\ell \geq 1}  \alpha^{(0)}_{\ell,X}\left\|  \|X_{t-\ell}\right\|\|_2\\ 
       &\leq  
       C+   C\big(\sum\limits_{\ell \geq 1}  \alpha^{(0)}_{\ell,Y} +\sum\limits_{\ell \geq 1}  \alpha^{(0)}_{\ell,X}\big) <\infty.
  \end{align*}
   Hence, (\textbf{A5}) is satisfied with $ \overline{\psi}_t = C+   C\big(\sum\limits_{\ell \geq 1}  \alpha^{(0)}_{\ell,Y} +\sum\limits_{\ell \geq 1}  \alpha^{(0)}_{\ell,X}\big)$. 
   \item To verify (\textbf{A6}), remark that for the model (\ref{Model_Expo}), $\eta'(\lambda) =\frac{1}{B'(B^{-1}(\lambda))}$. Moreover, since $B$ is strictly increasing in $\eta$, $B^{-1}$ is also strictly increasing in $\lambda$. Then, from  (\textbf{A0}), 
$
  \left|\eta'(\lambda)\right| \leq \frac{1}{B'(B^{-1}(\underline{c}))}
$, for some $\underline{c}>0$,  
which implies that $\eta'$ is bounded. Thus, the function $\eta$ satisfies the Lipschitz condition, which shows that (\textbf{A6}) holds.
Clearly, the second part of (\textbf{A6}) is satisfied.

\item Now, let us show that  (\textbf{A7}) is satisfied.  
From (\textbf{A0}), for all $\theta \in \Theta$, we have 
  \begin{align*}
\left| h_{\alpha}(\lambda_t(\theta))\right| 
&=
\frac{(1+\alpha)}{B'(B^{-1}(\lambda_t(\theta)))}\Big|\underset{y =0}{\overset{\infty}{\sum}} (y-\lambda_t(\theta))g (y|\eta(\lambda_t(\theta)))^{\alpha+1}  +(Y_t-\lambda_t(\theta))g (Y_t|\eta(\lambda_t(\theta)))^{\alpha}  \Big|\\
&\leq
\frac{(1+\alpha)}{B'(B^{-1}(\underline{c}))}\Big(\underset{y =0}{\overset{\infty}{\sum}} (y+\lambda_t(\theta))g (y|\eta(\lambda_t(\theta)))  +(Y_t+\lambda_t(\theta))\Big)\\
&\leq C\left( Y_t+3\lambda_t(\theta) \right).
  \end{align*}
  Therefore, 
   \begin{equation}\label{eq1_Ex1} 
  \left\|h_\alpha(\lambda_t(\theta)) \right\|^2_\Theta \leq C
  \| Y_t+3\lambda_t(\theta)\|^2_\Theta \leq C ( Y^2_t+6Y_t \|\lambda_t(\theta)\|_\Theta + 9 \|\lambda_t(\theta)\|^2_\Theta ).
  \end{equation}
   By applying the H\"{o}lder's inequality to the second term of the right hand side of (\ref{eq1_Ex1} ), we get
   \begin{align*}
  \left\|\left\|h_\alpha(\lambda_t(\theta)) \right\|^2_\Theta\right\|_2
  &\leq C
  \left(
   \left\|Y^2_t\right\|_2+6\big\|Y_t \cdot 
  \|\lambda_t(\theta)\|_\Theta \big\|_2 + 9  
  \left\|\|\lambda_t(\theta)\|^2_\Theta\right\|_2
  \right)\\
  &\leq C
  \left(
   \left\|Y_t\right\|^2_4+6\|Y_t\|_4 \cdot
  \|\|\lambda_t(\theta)\|_\Theta\|_4 + 9  
  \left\|\|\lambda_t(\theta)\|_\Theta\right\|^2_4
  \right) \\
  &\leq C
  \left(
   C+6C
  \underline{\lambda}^{1/4} + 9  
  \underline{\lambda}^{1/2}
  \right):= \overline{h}_\alpha<\infty ~(\text{from (\ref{moment}) with $\epsilon>3$ and (\textbf{A2})}).
  \end{align*}
  To complete the verification of (\textbf{A7}), we need to impose the following regularity condition on the function $B$ (see also Kim and Lee (2019)):
  
  \medskip 
   
	(\textbf{B0}): $\sup_{\theta \in \Theta} \sup_{0 \leq \delta \leq 1} \Big| \frac{B''(\eta_{t,\delta}(\theta))}{B'( \eta_{t,\delta}(\theta))^3}\Big| \leq K$, for some $K>0$. 
	
	\medskip 

\noindent 
Let $ \lambda_{t,\delta}(\theta)= \delta \lambda_t(\theta) + (1-\delta) \widehat{\lambda}_t(\theta)$ with $\theta \in \Theta$, $0\leq \delta \leq 1$; 
which implies that $\eta_{t,\delta}(\theta)=\eta(\lambda_{t,\delta}(\theta))$ is between $ \eta(\lambda_t(\theta))$ and $\eta(\widehat{\lambda}_t(\theta))$ since the function $\eta$ is monotone.
By using the condition (\textbf{B0}), we can proceed as in Kim and Lee (2019) to get for all $\theta \in \Theta$,
\begin{align*}
&\left| m_{\alpha}(\lambda_{t,\delta}(\theta))\right| 
\leq C Y^2_t+K Y_t+C B(\eta_{t,\delta}(\theta))^2+3K B(\eta_{t,\delta}(\theta))+C\\
&\hspace{0.8cm} \leq C Y^2_t+K Y_t+C \big(B(\eta_{t}(\theta)) +\left|B(\widehat \eta_{t}(\theta))- B(\eta_{t}(\theta))\right|\big)^2 +3K \big(B(\eta_{t}(\theta)) +\left|B(\widehat \eta_{t}(\theta))- B(\eta_{t}(\theta))\right|\big)+C\\
& \hspace{0.8cm} =C Y^2_t+K Y_t+C \big(\lambda_t(\theta) +|\widehat \lambda_t(\theta)-\lambda_t(\theta) |\big)^2+3K \big(\lambda_t(\theta) +|\widehat \lambda_t(\theta)-\lambda_t(\theta)|\big)+C.
\end{align*}
Thus, according to (\ref{App_A0}) and (\ref{App1_prime_A0}), for any $t \geq 1$, we have
\begin{multline*}
\left\| m_{\alpha}(\lambda_{t,\delta}(\theta))\right\|_\Theta \leq
 C Y^2_t+K Y_t+C \big(\left\| f_\theta(0)\right\|_\Theta+   2\sum\limits_{\ell \geq 1}  \alpha^{(0)}_{\ell,Y} Y_{t-\ell}+ 2\sum\limits_{\ell \geq 1}  \alpha^{(0)}_{\ell,X}\left\|  X_{t-\ell}\right\| \big)^2 \\
 +3K\big(\left\| f_\theta(0)\right\|_\Theta+   2\sum\limits_{\ell \geq 1}  \alpha^{(0)}_{\ell,Y} Y_{t-\ell}+ 2\sum\limits_{\ell \geq 1}  \alpha^{(0)}_{\ell,X}\left\|  X_{t-\ell}\right\| \big)+C.
\end{multline*}
Hence, from (\ref{moment}) with $\epsilon>3$, 
\begin{align*}
&\left\|\left\| m_{\alpha}(\lambda_{t,\delta}(\theta))\right\|_\Theta\right\|_2\\
&\hspace{1cm} \leq
 C \|Y^2_t\|_2+K \|Y_t\|_2+C \Big\|\big(\left\| f_\theta(0)\right\|_\Theta+   2\sum\limits_{\ell \geq 1}  \alpha^{(0)}_{\ell,Y} Y_{t-\ell}+ 2\sum\limits_{\ell \geq 1}  \alpha^{(0)}_{\ell,X}\left\|  X_{t-\ell}\right\|\big)^2\Big\|_2\\
 & \hspace{3cm} +3K\Big\|\big(\left\| f_\theta(0)\right\|_\Theta+   2\sum\limits_{\ell \geq 1}  \alpha^{(0)}_{\ell,Y} Y_{t-\ell}+ 2\sum\limits_{\ell \geq 1}  \alpha^{(0)}_{\ell,X}\left\|  X_{t-\ell}\right\|\big)\Big\|_2+C\\
 & \hspace{1cm}\leq
 C \|Y_t\|^2_4+K \|Y_t\|_2+C \Big\|\left\| f_\theta(0)\right\|_\Theta+   2\sum\limits_{\ell \geq 1}  \alpha^{(0)}_{\ell,Y} Y_{t-\ell}+ 2\sum\limits_{\ell \geq 1}  \alpha^{(0)}_{\ell,X}\left\|  X_{t-\ell}\right\|\Big\|^2_4\\
 & \hspace{3cm} +3K\big(\left\|\left\| f_\theta(0)\right\|_\Theta\right\|_2+   2\sum\limits_{\ell \geq 1}  \alpha^{(0)}_{\ell,Y}  \|Y_{t-\ell}\|_2+ 2\sum\limits_{\ell \geq 1}  \alpha^{(0)}_{\ell,X}\|\left\|  X_{t-\ell}\right\|\|_2\big)+C\\
  &\hspace{1cm} \leq
 C +CK +C \Big(\left\|\left\| f_\theta(0)\right\|_\Theta\right\|_4+   2\sum\limits_{\ell \geq 1}  \alpha^{(0)}_{\ell,Y}  \|Y_{t-\ell}\|_4+ 2\sum\limits_{\ell \geq 1}  \alpha^{(0)}_{\ell,X}\|\left\|  X_{t-\ell}\right\|\|_4\Big)^2\\
 & \hspace{3cm} +3K\big(C+   C \big(\sum\limits_{\ell \geq 1}  \alpha^{(0)}_{\ell,Y} + \sum\limits_{\ell \geq 1}  \alpha^{(0)}_{\ell,X}\big) \big)+C\\
 & \hspace{1cm} \leq
 C  +C \Big(C+   C \big(\sum\limits_{\ell \geq 1}  \alpha^{(0)}_{\ell,Y} + \sum\limits_{\ell \geq 1}  \alpha^{(0)}_{\ell,X}\big)\Big)^2
 +C\big(C+   C \big(\sum\limits_{\ell \geq 1}  \alpha^{(0)}_{\ell,Y} + \sum\limits_{\ell \geq 1}  \alpha^{(0)}_{\ell,X}\big)\big)+C:=\overline m_\alpha<\infty.
\end{align*}
Thus (\textbf{A7}) is satisfied.
  \end{itemize}
  

  \subsection{A particular case of linear models: INGARCH-X}\label{Sec_Exempl_INGARCH-X} 
 As second example, consider the Poisson-INGARCH-X model defined by
 	\begin{equation}\label{P_INGARCH_X}
            	Y_{t}|\mathcal{F}_{t-1} \sim \mathcal{P}(\lambda_{t})~~\text{with}~~
            	 \lambda_{t}=  	\alpha^*_{0} + \sum_{i=1}^{q^*} \alpha^*_{i} Y_{t-i} +  \sum_{i=1}^{p^*} \beta^*_{i} \lambda_{t-i}                                                       +v(X_{t-1}), 
    \end{equation}  
    where
   $ \alpha^*_{0}>0$, $\alpha^*_{1},\cdots, \alpha^*_{q^*} , \beta^*_{1},\cdots,\beta^*_{p^*} \geq 0$ and $v$ is a non-negative function defined on $\R^{d_x}$. 
 This class of models has already been studied within the maximum-likelihood framework; see for instance Agosto \textit{et al.} (2016) and Pedersen and Rahbek (2018). 
%
 Similarly, one can define the NB-INGARCH-X and BIN-INGARCH-X models with the negative binomial distribution and the Bernoulli distribution, respectively. 
 Without loss of generality, we assume that the components of the exogenous covariate vector are non-negative and that $v(x)$ is a linear function in $x$. 
 More precisely, there exists $(\gamma^*_1,\cdots,\gamma^*_{d_x}) \in [0,\infty)^{d_x}$  such that $v(x)=\sum\limits_{i = 1}^{d_x} \gamma^*_i x_{i}$, for any $x=(x_1,\cdots,x_{d_x}) \in [0,\infty)^{d_x}$. 
   The true parameter of the model is $\theta^*=( \alpha^*_{0},\alpha^*_{1},\cdots, \alpha^*_{q^*} , \beta^*_{1},\cdots,\beta^*_{p^*}, \gamma^*_1,\cdots,\gamma^*_{d_x})$; 
 therefore, $\Theta$ is a compact subset of $ (0,\infty) \times [0,\infty)^{p^* + q^*+d_x}$.  
 If $\sum\limits_{i=1}^{p^*} \beta_i < 1$, then we can find  two sequences 
$(\psi_k(\theta^*))_{k \geq 0}$ and $(\gamma_k(\theta^*))_{k \geq 1}$ such that

\[ \lambda_t = \psi_0(\theta^*) + \sum_{k \geq1} \psi_k(\theta^*)Y_{t-k}+ \sum_{k \geq1} \gamma'_k(\theta^*)X_{t-k}; \]

which implies that 
$f_{\theta^*}(y;x)=\psi_0(\theta^*) + \sum_{k \geq1} \psi_k(\theta^*)y_{k} + \sum_{k \geq1} \gamma'_k(\theta^*)x_{k}$, for any $(y,x) \in \N^{\N}_0 \times {(\R^{d_x})}^{\N}$. 
 As pointed out in Aknouche and Francq (2020),  by using the arguments of Theorem 1 in Francq and Thieu (2019), a sufficient condition on the model (\ref{P_INGARCH_X}) to be identifiable is that $Y_t$ is not a measurable function of $\{ X_s,\ s \in \Z\}$. 
%
Let us impose here a Markov-structure on the set of the covariates. Assume that $X_t=g(X_{t-1}, \varepsilon_t)$ for some function $g(x, \varepsilon_t)$ with values in $[0,\infty)^{d_x}$ and where ($\varepsilon_t$) is a sequence of  \textit{i.i.d} random variables. To ensure the stability of the exogenous covariates,  we make the following assumption (see also Agosto \textit{et al.} (2016)). 
    
    \medskip 
    
    (\textbf{B1}): $\E\left[\left\|g(x, \varepsilon_t)-g(x', \varepsilon_t) \right\|^s\right] \leq \rho \left\| x-x' \right\|^s$
      and 
      $\E\left[\left\|g(0, \varepsilon_t) \right\|^s\right] < \infty$, for some $(s,\rho) \in [1,\infty) \times (0,1)$ and for all $x, x' \in [0,\infty)^{d_x}$.
    
    \medskip 

\noindent 
To check the conditions of Theorems \ref{th1} and \ref{th2}, define the compact set $\Theta$ as follows:
    \begin{multline}\label{def_Theta}
\Theta=\Big\{
\theta=(\alpha_0,\alpha_1,\cdots,\alpha_{q^*},\beta_1,\cdots,\beta_{p^*},\gamma_1,\cdots,\gamma_{d_x}) \in [0,\infty)^{p^*+q^*+d_x+1} \big/\
0<\alpha_{L}\leq \alpha_{0} \leq \alpha_{U},
\\
 \sum\limits_{i = 1}^{q^*} \alpha_i+ \sum\limits_{i = 1}^{p^*}\beta_i< 1- \epsilon \text{ with } \epsilon=\max \big\{0, \big(1-\sum\limits_{i = 1}^{p^*}\beta_i\big)(\rho-\alpha_{1}) \big\}
~\text{ and } \max_{1 \leq i \leq d_x}\{ \gamma_i \} \leq \alpha_V 
 \Big\},
\end{multline}
for some $\alpha_L, \alpha_U,\alpha_V, \epsilon>0$.
\begin{itemize}
	\item  Under the assumption (\textbf{B1}), if $\theta^* \in \Theta$ (see the definition in \ref{def_Theta}), \textbf{A}$_i (\Theta)$ ($i=0,1,2$) holds; and
	Proposition \ref{exp_existence} establishes the existence of a  $\tau$-weakly dependent stationary and ergodic solution $Y^*_t=(Y_t,\lambda_t,X_t)$ of the model (\ref{P_INGARCH_X}) as well as the NB-INGARCH-X and BIN-INGARCH-X models, satisfying $\E\|Y^*_t\|<\infty$. In the case of the Poisson-INGARCH-X model, we refer to the second part of the proof of Theorem 1 in Agosto \textit{et al.} (2016) for the existence of the $s-$order moment with $s> 1$.  
Thus, according to (\textbf{B1}) (with $s> 3$ ), the condition (\ref{moment})  (with $\epsilon > 3$) holds,  and consequently (\textbf{A2}) holds for the Poisson-INGARCH-X model (see subsection \ref{Sec_Exempl_Ex}, and for sufficient conditions in a large class of models, including  NB-INGARCH-X and BIN-INGARCH-X).
In addition, $f_{\theta}(y; x)  \geq \alpha_{0} \geq \alpha_L :=\underline{c}$ for all $\theta \in \Theta$ and  $(y,x) \in \N^\N_0 \times {(\R^{d_x})}^\N$,  which shows that the second part of (\textbf{A0}) holds.  
	 Since the one-parameter exponential family includes the Poisson distribution, by using (\ref{moment}),  (\textbf{A0})  and \textbf{A}$_0 (\Theta)$, we can go along similar lines as in the general model (\ref{Model_Expo}) to show that the assumptions (\textbf{A3}),  (\textbf{A4}), (\textbf{A5}) and (\textbf{A6}) hold.   
  
\item 
Let us check the assumption (\textbf{A7}). If (\textbf{B0}) holds, one can go along similar lines as in (\ref{Model_Expo}) to get (\textbf{A7}).  Let us show that  (\textbf{B0}) holds. 
Remark that  model (\ref{P_INGARCH_X}) is  a particular case of model (\ref{Model_Expo}) where 
 $\lambda_t=B(\eta_t)=\rm{e}^{\eta_t}$. 
Thus, for any $\delta \in [0,1]$, $\theta \in \Theta$,  we have
\begin{align*}
\Big| \frac{B'(\eta_{t,\delta}(\theta))}{B''(\eta_{t,\delta}(\theta))^3} \Big| = \frac{1}{\lambda^2_{t,\delta}(\theta)} \leq \frac{1}{\alpha^2_{L}} \coloneqq K, 
\end{align*}
where $\lambda_{t,\delta}(\theta)$ is between $\widehat \lambda_{t}(\theta)$ and  $\lambda_{t}(\theta)$. Hence, (\textbf{B0}) is satisfied.
\end{itemize}


 \section{Simulation and results}
This section presents some simulation study to assess the efficiency and the robustness of the MDPDE.  
 We compare the performances of the MDPDE with those of the MLE ($\alpha =0$) for some dynamic models satisfying  (\ref{Model}). 
 To this end, the stability of the estimators under contaminated data will be studied.  
 For each model considered, the results are based on $100$ replications of Monte Carlo simulations of sample sizes $n=500,\,1000$.
  The sample mean and the mean square error (MSE) of the estimators will be applied as evaluation criteria. 



 %
  
\subsection{Poisson-INGARCH-X process}
 Consider the Poisson-INGARCH-X defined by 
 
 	\begin{equation}\label{Poisson_INGARCH_X}
            	Y_{t}|\mathcal{F}_{t-1} \sim \mathcal{P}(\lambda_{t})~~\text{with}~~
            	 \lambda_{t}=  	\alpha_{0} + \alpha_{1} Y_{t-1} +   \beta \lambda_{t-1}  
            	                                                      + \gamma |X_{t-1}|, 
    \end{equation}  
    where 
      $\theta^*=(\alpha_{0},\alpha_{1}, \beta,\gamma) \in \Theta \subset  (0,\infty) \times [0,\infty)^3$ is the true parameter  
    and 
    $X_t$ follows an ARCH$(1)$ process given by
    \[
    X_{t}|\mathcal{F}_{t-1} \sim \mathcal N(0,\sigma^2_t) ~~\text{with}~~\sigma^2_t=\omega_0+ \omega_1 X^2_{t-1} 
            ,~~\omega_0 >0,\,\omega_1 \geq 0.
    \] 
     For this model, we set $\theta^*=(0.10,0.15,0.80,0.03)$. This scenario is related and close to the real data example (see below). 
   We first consider the case where the data are not contaminated by outliers. We generate a trajectory of model (\ref{Poisson_INGARCH_X}) with $(\omega_0,\omega_1) =(1,0.5)$.   
   For different values of $\alpha$, the parameter estimates and their corresponding MSEs (shown in parentheses) are summarized in
   Table \ref{Res1.Poisson_INGARCH_X}. In the table, the minimal MSE for each component of  $\widehat{\theta}_{\alpha,n}$  is indicated by the symbol $^{\bullet}$ (some values of the MSE have been rounded up, e.g. in Table\ref{Res1.Poisson_INGARCH_X}, $n=500$, $\alpha=0.10$, the last column, in parentheses, the value is 0.088 instead of 0.09).  
    These results show that the MLE has the minimal MSEs for all the parameters, except for $\gamma$ when $n=500$ and $\alpha_1$ when $n=1000$.
    For the parameters $\gamma$ (when $n=500$) and $\alpha_1$ (when $n=1000$), the MDPDE with $\alpha=0.1$ has the minimal MSE, but the values are close to those of the MLE. 
        One can also observe that the MSEs of the MDPDEs increase with $\alpha$. 
         These findings confirm the fact that the  MLE generally outperforms the MDPDE when no outliers exist. 
   However, as $n$ increases, the performances of the MDPDE  increase for each $\alpha$.\\ 
   Now, we evaluate the robustness of the estimators by considering the case where the data are contaminated by additive outliers.  Assume that we observe the contaminated process $Y_{c,t}$ such that $Y_{c,t}=Y_{t}+P_t Y_{0,t}$, where $Y_{t}$  is generated from (\ref{Poisson_INGARCH_X}), $P_t$ is an \textit{i.i.d} Bernoulli random variable with a success probability $p$ and $Y_{0,t}$ is an \textit{i.i.d} Poisson random variable with a mean $\mu$. In the sequel, the variables $P_t$, $Y_{0,t}$ and $Y_{t}$ are assumed to be independent.  For $p=0.02$ and $\mu=10$, the corresponding results are summarized in Table \ref{Res2.Poisson_INGARCH_X}. 
 From this table, one can see that the MDPDE has smaller MSEs than the MLE, except for the estimations obtained 
 with $\alpha=0.75$ and $\alpha=1$ (see for instance $\widehat \beta$); 
 which indicates that the MDPDE is more robust to outliers and overall outperforms the MLE in such cases.  
  We also observe that the selected optimal value of $\alpha$ decreases as $n$ increases for all parameters. 
   
 \begin{table}[h!]
\scriptsize
\centering
\caption{\it Sample mean and MSE$\times 10^2$ of the estimators for the Poisson-INGARCH-X model (\ref{Poisson_INGARCH_X}) with $\theta^*=(0.1,0.15,0.8, 0.03)$:  the case without outliers in the data.}
\label{Res1.Poisson_INGARCH_X}
\vspace{.2cm}
\hspace*{-.3cm}
\begin{tabular}{c llll  c llll}
\Xhline{.9pt}
\rule[0cm]{0cm}{.4cm}
 &
  \multicolumn{4} {l} {$n=500$}&&
  \multicolumn{4} {l} {$n=1000$}  \\
\cline{2-5}\cline{7-10}
\rule[0cm]{0cm}{.35cm}
  $\alpha$
  &$\widehat \alpha_{0}$&$\widehat \alpha_{1}$&$\widehat \beta$& $\widehat \gamma$&&
  $\widehat \alpha_{0}$&$\widehat \alpha_{1}$&$\widehat \beta$&$\widehat \gamma$\\
\Xhline{.72pt}
\rule[0cm]{0cm}{.3cm}   
                        $0$
                        &$0.137(1.86)^{\bullet}$&$0.152(0.10)^{\bullet}$&$0.780(0.60)^{\bullet}$&$0.035(0.09)$&$\bf \cdot$&
                         $0.120(0.39)^{\bullet}$&$0.148(0.06)$&$0.793(0.16)^{\bullet}$&$0.034(0.06)^{\bullet}$\\ 
 \rule[0cm]{0cm}{0cm} $0.10$  
                          &$0.139(2.05)$&$0.152(0.11)$&$0.779(0.67)$&$0.035(0.09)^{\bullet}$&$\bf \vdots$&
                           $ 0.121(0.39)$&$0.148(0.06)^{\bullet}$&$0.793(0.16)$&$0.034(0.06)$\\ 
 \rule[0cm]{0cm}{0cm}  $0.20$ 
                        &$0.140(2.25)$&$0.152(0.11)$&$0.779(0.74)$&$0.034(0.09)$&$\bf \vdots$&
                         $0.122(0.41)$&$0.148(0.07)$&$0.792(0.16)$&$0.034(0.06)$\\
 \rule[0cm]{0cm}{.0cm}  $0.30$  
                        &$0.143(2.46)$&$0.152(0.12)$&$0.778(0.81)$&$0.034(0.09)$&$\bf \vdots$&
                         $0.123(0.43)$&$0.149(0.07)$&$0.791(0.17)$&$0.034(0.06)$\\    
     
 \rule[0cm]{0cm}{.0cm}  $0.40$  
                        &$0.144(2.64)$&$0.152(0.12)$&$0.777(0.87)$&$0.034(0.10)$&$\bf \vdots$&
                         $0.123(0.45)$&$0.149(0.07)$&$0.791(0.18)$&$0.034(0.06)$\\    
 \rule[0cm]{0cm}{.0cm}  $0.50$  
                        &$0.146(2.81)$&$0.152(0.13)$&$0.776(0.93)$&$0.034(0.10)$&$\bf \vdots$&
                         $0.124(0.47)$&$0.149(0.07)$&$0.790(0.19)$&$0.034(0.07)$\\  
 \rule[0cm]{0cm}{.4cm}  $0.75$  
                        &$0.150(3.09)$&$0.152(0.14)$&$0.775(1.03)$&$0.033(0.10)$&$\bf \vdots$&
                         $0.127(0.53)$&$0.150(0.09)$&$0.788(0.22)$&$0.034(0.07)$\\  
 \rule[0cm]{0cm}{.4cm}  $1.00$  
                        &$0.153(3.23)$&$0.152(0.16)$&$0.774(1.10)$&$0.032(0.10)$&$\bf \cdot$&
                         $0.129(0.59)$&$0.151(0.10)$&$0.787(0.26)$&$0.034(0.08)$\\                 
       \Xhline{.9pt}
\end{tabular}
\end{table}
 
 
 \begin{table}[h!]
\scriptsize
\centering
\caption{\it Sample mean and MSE$\times 10^2$ of the estimators for the Poisson-INGARCH-X model (\ref{Poisson_INGARCH_X}) with $\theta^*=(0.1,0.15,0.8, 0.03)$: the case in which the data are contaminated by outliers.}
\label{Res2.Poisson_INGARCH_X}
\vspace{.2cm}
\hspace*{-.3cm}
\begin{tabular}{c llll  c llll}

\Xhline{.9pt}
\rule[0cm]{0cm}{.4cm}
 &
  \multicolumn{4} {l} {$n=500$}&&
  \multicolumn{4} {l} {$n=1000$}  \\
\cline{2-5}\cline{7-10}
\rule[0cm]{0cm}{.35cm}
  $\alpha$
  &$\widehat \alpha_{0}$&$\widehat \alpha_{1}$&$\widehat \beta$& $\widehat \gamma$&&
  $\widehat \alpha_{0}$&$\widehat \alpha_{1}$&$\widehat \beta$&$\widehat \gamma$\\
\Xhline{.72pt}
\rule[0cm]{0cm}{.3cm}   
                        $0$
                        &$0.279(26.90)$&$0.105(0.34)$&$0.775(3.53)$&$0.045(0.30)$&$\bf \cdot$&
                        $0.140(0.75)$&$0.107(0.25)$&$0.827(0.27)$&$0.044(0.12)$
                         \\ 
 \rule[0cm]{0cm}{.4cm} $0.10$  
                          &$0.186(14.81)$&$0.105(0.30)$&$0.805(1.88)$&$0.041(0.18)$&$\bf \vdots$&
                          $0.109(0.49)^{\bullet}$&$0.109(0.22)$&$0.829(0.23)^{\bullet}$&$0.043(0.09)^{\bullet}$
                          \\  
 \rule[0cm]{0cm}{.4cm}  $0.20$ 
                        & $0.200(21.91)$&$0.105(0.30)$&$0.799(2.73)$&$0.041(0.16)$&$\bf \vdots$&
                        $0.104(0.52)$&$0.110(0.21)^{\bullet}$&$0.828(0.24)$&$0.043(0.09)$
                          \\
 \rule[0cm]{0cm}{.4cm}  $0.30$  
                        &$0.160(10.67)^{\bullet}$&$0.104(0.31)$&$0.812(1.58)^{\bullet}$&$0.041(0.16)^{\bullet}$&$\bf \vdots$&
                        $0.103(0.56)$&$0.110(0.22)$&$0.827(0.25)$&$0.043(0.10)$
                         \\    
     
 \rule[0cm]{0cm}{.4cm}  $0.40$  
                        &$0.185(17.23)$&$0.104(0.31)$&$0.802(2.27)$&$0.041(0.16)$&$\bf \vdots$&
                        $0.104(0.58)$&$0.110(0.22)$&$0.826(0.26)$&$0.043(0.10)$
                          \\  
 \rule[0cm]{0cm}{.4cm}  $0.50$  
                        &$0.187(17.17)$&$0.105(0.31)$&$0.801(2.27)$&$0.041(0.16)$&$\bf \vdots$&
                        $0.106(0.63)$&$0.110(0.22)$&$0.825(0.27)$&$0.043(0.10)$
                        
                          \\ 
 \rule[0cm]{0cm}{.4cm}  $0.75$  
                        &$0.200(17.94)$&$0.107(0.30)$&$0.792(2.72)$&$0.041(0.17)$&$\bf \vdots$&
                         $0.109(0.69)$&$0.111(0.22)$&$0.822(0.29)$&$0.0441(0.11)$
                          \\  
 \rule[0cm]{0cm}{.4cm}  $1.00$  
                        &$0.238(26.31)$&$0.109(0.30)^{\bullet}$&$0.775(3.55)$&$0.043(0.19)$&$\bf \cdot$&
                        $0.112(0.78)$&$0.113(0.22)$&$0.818(0.32)$&$0.045(0.12)$
                         \\                  
       \Xhline{.9pt}
\end{tabular}
\end{table}



\subsection{NB-INGARCH-X process}
 Consider the negative-binomial-INGARCH-X (NB-INGARCH-X) model defined by 
 	\begin{equation}\label{NB_INGARCH_X}
            	Y_{t}|\mathcal{F}_{t-1} \sim NB(r,p_{t})~\text{with}~
            	 r\frac{(1-p_t)}{p_t}=\lambda_{t}=  	\alpha_{0} + \alpha_{1} Y_{t-1} +  \beta \lambda_{t-1}                                                       + \gamma_1 \exp(X_{1,t-1}) + \gamma_2 \textrm{\dsrom{1}}_{\{X_{2,t-1}<0\}} |X_{2,t-1}|, 
    \end{equation}  
    where 
     $\theta^*=(\alpha_{0},\alpha_{1}, \beta,\gamma_1,\gamma_2)  \in  \Theta \subset (0,\infty) \times [0,\infty)^4$ is the true parameter, $(X_{1,t},X_{2,t}) $ is the covariate vector, 
     $\{ X_{i,t},\, t\geq 1 \}$ (for $i=1,2$) is an autoregressive process satisfying   
   \[
      X_{i,t}=\varphi_i X_{i,t-1} + \varepsilon_{i,t}~ (\text{with } 0<\varphi_i <1 \text{ and }  \varepsilon_{i,t} \text{ is a Gaussian white noise}),
   \] 
  $\textrm{\dsrom{1}}_{\{\cdot\}}$ denotes the indicator function  
     and
     $NB(r,p)$ denotes the negative binomial distribution with parameters $r$ and $p$. 
For this model, we consider the cases of $(\varphi_1,\varphi_2)=(1/3,1/2)$, $r=8$ and $\theta^*=(0.5,0.2,0.4,0.1,0.3)$. 
We first generate a data $Y_t$ of model (\ref{NB_INGARCH_X}) (without outliers). 
To evaluate the robustness of the estimators, we consider the contaminated data $Y_{c,t}$ (presence of outliers) as follows:  $Y_{c,t}=Y_{t}+P_t Y_{0,t}$, where $P_t$ is an \textit{i.i.d} Bernoulli random variable with a success probability $p=0.02$ and $Y_{0,t}$ is an \textit{i.i.d} $NB(5,0.4)$. 
The results are presented in Tables \ref{Res1.NB_INGARCH_X} and \ref{Res2.NB_INGARCH_X}.   
Once again, in the absence of outliers (see Table \ref{Res1.NB_INGARCH_X}), the MLE outperforms the MDPDE  and the efficiency of the MDPDE decreases as $\alpha$ increases. 
 When the data are contaminated by outliers (see Table \ref{Res2.NB_INGARCH_X}), one can see that the MDPDE has smaller MSEs than the MLE; that is, the MDPDE is more robust than the MLE.  
 As in the model (\ref{Poisson_INGARCH_X}), when $n$ increases, the symbol $^{\bullet}$ overall tends to move upwards.

 \begin{table}[h!]
\scriptsize
\centering
\caption{\it Sample mean and MSE$\times 10^2$ of the estimators for the NB-INGARCH-X model (\ref{NB_INGARCH_X})  with $\theta^*=(0.5,0.2, 0.4,0.1,0.3)$: the case without outliers in the data.}
\label{Res1.NB_INGARCH_X}
\vspace{.2cm}
{\setlength{\tabcolsep}{2.pt}
\begin{tabular}{c lllll  c lllll}
\Xhline{.9pt}
\rule[0cm]{0cm}{.4cm}
 &
  \multicolumn{5} {l} {$n=500$}&&
  \multicolumn{5} {l} {$n=1000$}  \\
\cline{2-6}\cline{8-12}
\rule[0cm]{0cm}{.35cm}
  $\alpha$
  &$\widehat \alpha_{0}$&$\widehat \alpha_{1}$&$\widehat \beta$& $\widehat \gamma_1$& $\widehat \gamma_2$&&
  $\widehat \alpha_{0}$&$\widehat \alpha_{1}$&$\widehat \beta$&$\widehat \gamma_1$& $\widehat \gamma_2$\\
\Xhline{.72pt}
\rule[0cm]{0cm}{.3cm}   
                        $0$                        &$0.508(3.64)^{\bullet}$&$0.200(0.21)^{\bullet}$&$0.395(1.55)^{\bullet}$&$0.102(0.09)^{\bullet}$&$0.308(0.90)^{\bullet}$&$\bf \cdot$&
                           $0.507(1.79)^{\bullet}$&$0.192(0.13)^{\bullet}$&$0.405(0.71)^{\bullet}$&$0.102(0.06)^{\bullet}$&$0.300(0.51)^{\bullet}$\\ 
 \rule[0cm]{0cm}{.4cm} $0.10$  
                          &$0.509(3.75)$&$0.200(0.22)$&$0.394(1.60)$&$0.103(0.09)$&$0.309(0.91)$&$\bf \vdots$&
                           $0.507(1.92)$&$0.192(0.14)$&$0.405(0.76)$&$0.101(0.06)$&$0.299(0.52)$\\
 \rule[0cm]{0cm}{.4cm}  $0.20$ 
                         &$0.510(3.93)$&$0.201(0.23)$&$0.393(1.68)$&$0.103(0.10)$&$0.309(0.94)$&$\bf \vdots$&
                          $0.507(2.05)$&$0.193(0.14)$&$0.405(0.82)$&$0.101(0.06)$&$0.298(0.54)$\\
 \rule[0cm]{0cm}{.4cm}  $0.30$  
                        &$0.511(4.17)$&$0.201(0.24)$&$0.391(1.78)$&$0.104(0.10)$&$0.310(0.98)$&$\bf \vdots$&
                         $0.507(2.19)$&$0.193(0.15)$&$0.406(0.88)$&$0.101(0.06)$&$0.298(0.56)$\\  
     
 \rule[0cm]{0cm}{.4cm}  $0.40$  
                         &$0.512(4.45)$&$0.201(0.25)$&$0.390(1.89)$&$0.104(0.11)$&$0.310(1.04)$&$\bf \vdots$&
                          $0.507(2.31)$&$0.193(0.15)$&$0.406(0.95)$&$0.101(0.07)$&$0.298(0.59)$\\    
 \rule[0cm]{0cm}{.4cm}  $0.50$  
                         &$0.512(4.69)$&$0.201(0.26)$&$0.389(1.99)$&$0.105(0.12)$&$0.312(1.11)$&$\bf \vdots$&
                          $0.507(2.43)$&$0.193(0.16)$&$0.406(0.99)$&$0.100(0.07)$&$0.298(0.62)$\\ 
 \rule[0cm]{0cm}{.4cm}  $0.75$  
                         &$0.512(5.29)$&$0.202(0.30)$&$0.387(2.24)$&$0.106(0.14)$&$0.312(1.31)$&$\bf \vdots$&
                          $0.506(2.67)$&$0.192(0.18)$&$0.407(1.10)$&$0.100(0.08)$&$0.298(0.69)$\\
 \rule[0cm]{0cm}{.4cm}  $1.00$  
                         &$0.514(5.97)$&$0.203(0.35)$&$0.384(2.52)$&$0.107(0.16)$&$0.315(1.54)$&$\bf \vdots$&
                          $0.505(2.86)$&$0.193(0.19)$&$0.408(1.20)$&$0.100(0.08)$&$0.299(0.77)$\\                 
       \Xhline{.9pt}
\end{tabular}
}
\end{table}


 \begin{table}[h!]
\scriptsize
\centering
\caption{\it Sample mean MSE$\times 10^2$ of the estimators for the NB-INGARCH-X model (\ref{NB_INGARCH_X}) with $\theta^*=(0.5,0.2, 0.3,0.1,0.3)$:  the case in which the data are contaminated by outliers.}
\label{Res2.NB_INGARCH_X}
\vspace{.2cm}
{\setlength{\tabcolsep}{2pt}
\begin{tabular}{c lllll  c lllll}

\Xhline{.9pt}
\rule[0cm]{0cm}{.4cm}
 &
  \multicolumn{5} {l} {$n=500$}&&
  \multicolumn{5} {l} {$n=1000$}  \\
\cline{2-6}\cline{8-12}
\rule[0cm]{0cm}{.35cm}
  $\alpha$
  &$\widehat \alpha_{0}$&$\widehat \alpha_{1}$&$\widehat \beta$& $\widehat \gamma_1$& $\widehat \gamma_2$&&
  $\widehat \alpha_{0}$&$\widehat \alpha_{1}$&$\widehat \beta$&$\widehat \gamma_1$& $\widehat \gamma_2$\\
\Xhline{.72pt}
\rule[0cm]{0cm}{.3cm}   
                        $0$
                          &$0.608(10.57)$&$0.140(0.66)$&$0.433(2.91)$&$0.107(0.20)$&$0.303(1.02)$&$\bf \cdot$&
                            $0.619(5.35)$&$0.153(0.39)$&$0.416(1.35)$&$0.103(0.09)$&$0.312(0.65)$\\ 
 \rule[0cm]{0cm}{.4cm} $0.10$  
                          &$0.546(5.97)$&$0.146(0.51)$&$0.433(2.19)$&$0.107(0.16)$&$0.310(0.79)$&$\bf \vdots$&
                            $0.564(2.70)$&$0.160(0.29)$&$0.414(0.92)$&$0.104(0.06)^{\bullet}$&$0.317(0.60)^{\bullet}$\\
 \rule[0cm]{0cm}{.4cm}  $0.20$ 
                         &$0.526(5.25)$&$0.148(0.49)$&$0.433(2.09)^{\bullet}$&$0.107(0.15)^{\bullet}$&$0.311(0.77)^{\bullet}$&$\bf \vdots$&
                          $0.546(2.29)$&$0.161(0.28)^{\bullet}$&$0.413(0.89)^{\bullet}$&$0.103(0.06)$&$0.318(0.61)$\\
 \rule[0cm]{0cm}{.4cm}  $0.30$  
                        &$0.518(5.22)^{\bullet}$&$0.148(0.48)$&$0.432(2.14)$&$0.106(0.15)$&$0.311(0.81)$&$\bf \vdots$&
                          $0.540(2.21)^{\bullet}$&$0.162(0.28)$&$0.412(0.92)$&$0.103(0.06)$&$0.317(0.63)$\\  
     
 \rule[0cm]{0cm}{.4cm}  $0.40$  
                         &$0.518(5.37)$&$0.149(0.48)^{\bullet}$&$0.429(2.19)$&$0.106(0.16)$&$0.312(0.86)$&$\bf \vdots$&
                            $0.537(2.25)$&$0.162(0.29)$&$0.411(0.97)$&$0.103(0.07)$&$0.316(0.65)$\\    
 \rule[0cm]{0cm}{.4cm}  $0.50$  
                         &$0.518(5.62)$&$0.150(0.49)$&$0.428(2.31)$&$0.106(0.16)$&$0.311(0.91)$&$\bf \vdots$&
                          $0.537(2.35)$&$0.163(0.29)$&$0.410(1.03)$&$0.103(0.07)$&$0.316(0.66)$\\ 
 \rule[0cm]{0cm}{.4cm}  $0.75$  
                         &$0.532(6.63)$&$0.153(0.50)$&$0.415(2.71)$&$0.106(0.19)$&$0.315(1.10)$&$\bf \vdots$&
                          $0.538(2.63)$&$0.165(0.29)$&$0.406(1.17)$&$0.103(0.08)$&$0.315(0.71)$\\
 \rule[0cm]{0cm}{.4cm}  $1.00$  
                         &$0.535(7.05)$&$0.157(0.51)$&$0.409(2.93)$&$0.106(0.21)$&$0.318(1.26)$&$\bf \vdots$&
                          $0.540(2.88)$&$0.168(0.29)$&$0.402(1.29)$&$0.104(0.09)$&$0.315(0.75)$\\                 
       \Xhline{.9pt}
\end{tabular}
}
\end{table}




 \subsection{A $1$-knot  dynamic model}

 We consider the $1$-knot nonlinear dynamic model defined by (see also Davis and Liu (2016))  
\begin{equation}\label{model_noeud}
    Y_{t}|\mathcal{F}_{t-1} \sim \mathcal P (\lambda_{t})~~\text{with}~~  \lambda_{t} =
    \alpha_{0} + \alpha_{1}Y_{t-1} + \alpha_{2}\lambda_{t-1} + \beta(Y_{t-1}-\xi^*)^+ ,
   \end{equation}
 where $\alpha_{0}>0$, $\alpha_{1},\alpha_{2}, \beta\geq 0$,   
 $\xi^*$ is a non-negative integer (so-called knot) and $x^+ = \max(x,0)$ is the positive part of $x$.  
 The model (\ref{model_noeud}) is a particular case of the models (\ref{Model}) and (\ref{Model_Expo})  with $X_t \equiv constant$. 
 The true parameter is $\theta^*=(\alpha_{0}, \alpha_{1},\alpha_{2}, \beta)$. 
In this model, we consider the cases where $\xi^*=4$ and $\theta^*=(1,0.3,0.2,0.4)$. 
We generate a data $Y_t$ from (\ref{model_noeud}) and a contaminated data $Y_{c,t}$ such that $Y_{c,t}=Y_{t}+P_t Y_{0,t}$, where $P_t$ is an \textit{i.i.d} Bernoulli random variable with a success probability $p=0.006$ and $Y_{0,t}$ is an \textit{i.i.d} Poisson random variable with a mean $\mu=11$. 
 For each $\alpha$, the knot $\xi^*$ is estimated by minimizing the function $H_{\alpha,n}(\cdot)$ over the set of integer values $\{ 1, \cdots , \xi_{\max} \}$ where $\xi_{\max}$ is an upper bound of the true knot $\xi^*$ given by $\xi_{\max} = \max(Y_1, \cdots, Y_n)$. 
The estimation can be summarized as follows:
\begin{itemize}
	\item  For each $\xi \in \{1, \cdots , \xi_{\max} \}$ fixed, compute the MDPDE
of $\theta^*$ denoted  $\widehat{\theta}_{\alpha,n,\xi} $.
       
	\item  Estimate the knot by the relation:  
	$ 
	\widehat{\xi}_{\alpha,n} = \underset{ \xi \in \{ 1, \cdots , \xi_{\max} \}}{ \text{argmin}} H_{\alpha,n}(\widehat{\theta}_{\alpha,n,\xi}).
$
\end{itemize}
Some empirical statistics of the estimator $\widehat{\xi}_{\alpha,n}$ are reported in Table \ref{Stat_Estim_Xi}. 
These results show that the estimation of the knot is reasonably good in terms of the mean and the quantiles. In addition, the empirical probability of selecting the true knot increases with $n$.  
From Table \ref{Res1.noeud}, we can see that when the data are without outliers, the MLE displays the minimal MSE (except for the estimation of $\beta$);
%
%
whereas in the presence of outliers (see Table \ref{Res2.noeud}), the MDPDE outperforms the MLE. 
 Further, for the parameter $\beta$, the MDPDE with $\alpha = 1$ has the minimal MSE; these results reveal that the estimation of $\beta$ is more damaged than that of the other parameters.  
This can be explained by the fact that the term $(Y_{t-1}-\xi^*)^+$ (in the relation (\ref{model_noeud})) is very  sensitive to outliers.

\begin{table}[h!]
\scriptsize
\centering
\caption{ Some elementary statistics of the estimator $\widehat{\xi}_{\alpha,n}$ for the model (\ref{model_noeud}) without outliers.}
\label{Stat_Estim_Xi}
\vspace{.2cm}
\begin{tabular}{ccccccccc}
\Xhline{.9pt}
 & &&&&&&&\\
Sample size&Mean&SD&Min&$Q_1$&Med&$Q_3$&Max&$\mathbb P(\widehat{\xi}_{\alpha,n}=\xi^*)$\\
\Xhline{.6pt}
 \rule[0cm]{0cm}{.3cm} $500$
&$3.71$&$1.19$&$1.00$&$3.00$&$4.00$&$4.00$&$7.00$&$0.45$\\
\rule[0cm]{0cm}{.3cm}
 $1000$
&$3.96$&$0.78$&$2.00$&$4.00$&$4.00$&$4.00$&$6.00$&$0.54$\\

 \Xhline{.9pt}
\end{tabular}
\end{table} 
 

\begin{table}[h!]
\scriptsize
\centering
\caption{\it Sample mean and MSE$\times 10^2$ of the estimators for the $1-$knot dynamic model (\ref{model_noeud}) with $\theta^*=(1,0.3,0.2,0.4)$: the case without outliers in the data.}
\label{Res1.noeud}
\vspace{.2cm}
\hspace*{-.3cm}
\begin{tabular}{c llll  c llll}

\Xhline{.9pt}
\rule[0cm]{0cm}{.4cm}
 &
  \multicolumn{4} {l} {$n=500$}&&
  \multicolumn{4} {l} {$n=1000$}  \\
\cline{2-5}\cline{7-10}
\rule[0cm]{0cm}{.35cm}
  $\alpha$
  &$\widehat \alpha_{0}$&$\widehat \alpha_{1}$&$\widehat \alpha_{2}$& $\widehat \beta$&&
  $\widehat \alpha_{0}$&$\widehat \alpha_{1}$&$\widehat \alpha_{2}$&$\widehat \beta$\\
\Xhline{.72pt}
\rule[0cm]{0cm}{.3cm}   
                        $0$
                          &$1.052(3.96)^{\bullet}$&$0.288(1.02)^{\bullet}$&$0.172(0.88)^{\bullet}$&$0.405(1.96)$&$\bf \cdot$&
                           $1.028(2.38)^{\bullet}$&$0.290(0.29)^{\bullet}$&$0.189(0.60)^{\bullet}$&$0.382(1.33)$\\ 
 \rule[0cm]{0cm}{.4cm} $0.10$  
                          &$1.052(4.09)$&$0.283(1.09)$&$0.174(0.89)$&$0.404(1.84)^{\bullet}$&$\bf \vdots$&
                          $1.03(2.46)$&$0.289(0.30)$&$0.189(0.62)$&$0.387(1.29)$\\
 \rule[0cm]{0cm}{.4cm}  $0.20$ 
                         &$1.051(4.08)$&$0.280(1.13)$&$0.176(0.92)$&$0.409(1.92)$&$\bf \vdots$&
                          $1.034(2.53)$&$0.287(0.33)$&$0.188(0.64)$&$0.389(1.29)^{\bullet}$\\
 \rule[0cm]{0cm}{.4cm}  $0.30$  
                        &$1.048(4.16)$&$0.281(1.08)$&$0.178(0.95)$&$0.410(2.04)$&$\bf \vdots$&
                         $1.036(2.61)$&$0.285(0.34)$&$0.189(0.67)$&$0.388(1.38)$\\  
     
 \rule[0cm]{0cm}{.4cm}  $0.40$  
                         &$1.047(4.35)$&$0.280(1.13)$&$0.178(1.01)$&$0.414(2.09)$&$\bf \vdots$&
                          $1.037(2.75)$&$0.283(0.37)$&$0.189(0.69)$&$0.386(1.46)$\\    
 \rule[0cm]{0cm}{.4cm}  $0.50$  
                         &$1.048(4.59)$&$0.276(1.25)$&$0.179(1.07)$&$0.418(2.11)$&$\bf \vdots$&
                          $1.037(2.88)$&$0.284(0.37)$&$0.189(0.73)$&$0.390(1.50)$\\ 
 \rule[0cm]{0cm}{.4cm}  $0.75$  
                         &$1.045(5.08)$&$0.271(1.37)$&$0.181(1.25)$&$0.419(2.33)$&$\bf \vdots$&
                          $1.038(3.18)$&$0.282(0.45)$&$0.190(0.83)$&$0.397(1.66)$\\
 \rule[0cm]{0cm}{.4cm}  $1.00$  
                         &$1.046(5.88)$&$0.264(1.58)$&$0.183(1.44)$&$0.425(2.42)$&$\bf \vdots$&
                          $1.040(3.45)$&$0.275(0.56)$&$0.191(0.94)$&$0.393(1.93)$\\                 
       \Xhline{.9pt}
\end{tabular}
\end{table}


\begin{table}[h!]
\scriptsize
\centering
\caption{\it Sample mean and MSE$\times 10^2$ of the estimators for the $1-$knot dynamic model (\ref{model_noeud}) with $\theta^*=(1,0.3,0.2,0.4)$:  the case in which the data are contaminated by outliers.}
\label{Res2.noeud}
\vspace{.2cm}
\hspace*{-.3cm}
\begin{tabular}{c llll  c llll}

\Xhline{.9pt}
\rule[0cm]{0cm}{.4cm}
 &
  \multicolumn{4} {l} {$n=500$}&&
  \multicolumn{4} {l} {$n=1000$}  \\
\cline{2-5}\cline{7-10}
\rule[0cm]{0cm}{.35cm}
  $\alpha$
  &$\widehat \alpha_{0}$&$\widehat \alpha_{1}$&$\widehat \alpha_{2}$& $\widehat \beta$&&
  $\widehat \alpha_{0}$&$\widehat \alpha_{1}$&$\widehat \alpha_{2}$&$\widehat \beta$\\
\Xhline{.72pt}
\rule[0cm]{0cm}{.3cm}   
                        $0$
                          &$1.018(7.52)$&$0.258(1.19)$&$0.233(2.64)$&$0.134(9.78)$&$\bf \cdot$&
                           $1.020(2.77)$&$0.259(0.90)$&$0.233(0.88)$&$0.123(9.24)$\\ 
 \rule[0cm]{0cm}{.4cm} $0.10$  
                          &$0.998(5.22)^{\bullet}$&$0.269(0.93)$&$0.217(1.95)$&$0.137(9.64)$&$\bf \vdots$&
                           $1.004(2.39)^{\bullet}$&$0.261(0.74)$&$0.219(0.65)^{\bullet}$&$0.132(8.86)$\\
 \rule[0cm]{0cm}{.4cm}  $0.20$ 
                         &$1.011(5.32)$&$0.271(0.92)^{\bullet}$&$0.205(1.90)^{\bullet}$&$0.158(8.90)$&$\bf \vdots$&
                          $1.015(2.62$&$0.264(0.70)$&$0.208(0.65)$&$0.150(8.31)$\\
 \rule[0cm]{0cm}{.4cm}  $0.30$  
                        &$1.031(6.07)$&$0.270(1.01)$&$0.193(2.03)$&$0.187(8.28)$&$\bf \vdots$&
                          $1.028(2.84)$&$0.270(0.52)$&$0.198(0.70)$&$0.167(7.80)$\\  
     
 \rule[0cm]{0cm}{.4cm}  $0.40$  
                         &$1.048(6.72)$&$0.269(1.10)$&$0.183(2.14)$&$0.217(7.77)$&$\bf \vdots$&
                          $1.046(3.24)$&$0.269(0.53)$&$0.188(0.76)$&$0.190(7.18)$\\    
 \rule[0cm]{0cm}{.4cm}  $0.50$  
                         &$1.058(7.28)$&$0.267(1.14)$&$0.178(2.24)$&$0.234(7.40)$&$\bf \vdots$&
                          $1.059(3.68)$&$0.272(0.51)$&$0.179(0.85)$&$0.210(6.79)$\\ 
 \rule[0cm]{0cm}{.4cm}  $0.75$  
                         &$1.083(8.44)$&$0.268(1.22)$&$0.165(2.43)$&$0.289(6.53)$&$\bf \vdots$&
                          $1.080(4.17)$&$0.276(0.50)^{\bullet}$&$0.166(1.01)$&$0.262(6.04)$\\
 \rule[0cm]{0cm}{.4cm}  $1.00$  
                         &$1.089(9.14)$&$0.266(1.40)$&$0.162(2.62)$&$0.325(5.86)^{\bullet}$&$\bf \vdots$&
                          $1.102(4.98)$&$0.275(0.54)$&$0.154(1.19)$&$0.311(5.10)^{\bullet}$\\                 
       \Xhline{.9pt}
\end{tabular}
\end{table}



\section{Real data application}
\begin{figure}[h!]
\begin{center}
\includegraphics[height=10cm, width=13cm]{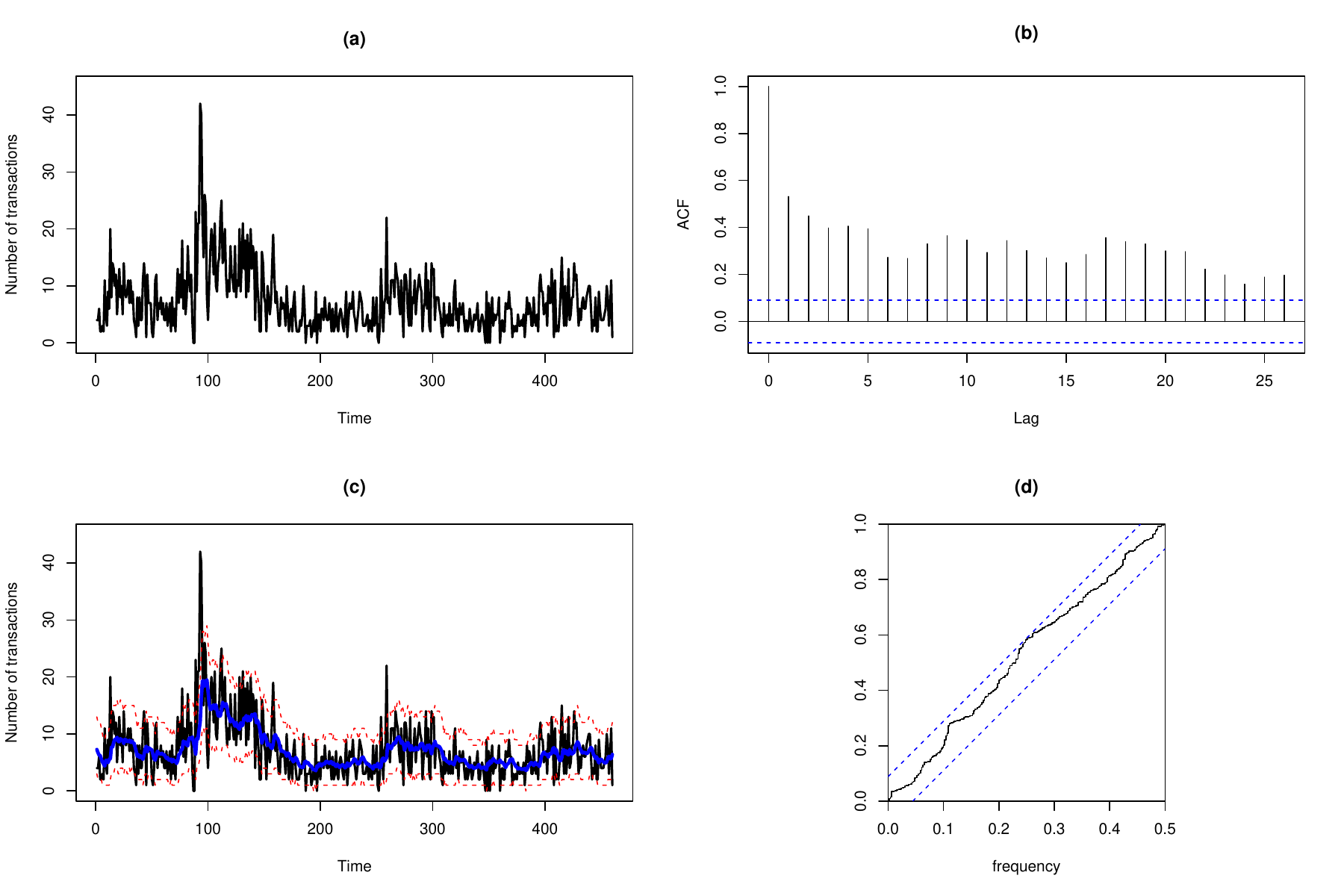}
\end{center}
\vspace{-1cm}
\caption{\it (a) Number of transactions per minute for the stock Ericsson B  during July 21, 2002. (b) Autocorrelation function of the transaction data. (c) Predicted  number of transactions per minute and the corresponding confidence bands at the $95\%$ nominal level (dotted curves) based on the relation (\ref{Estim_Appl}). (c) Cumulative periodogram plot of the Pearson residuals from the model (\ref{P_INGARCH_X(1,1)}).}
\label{Graphe}
\end{figure}

 The aim of this section is to apply the model (\ref{Model}) to analyze the number of transactions per minute for the stock Ericsson B during July 21, 2002. There are $460$ observations which represent trading from 09:35 to 17:14.
 This time series is a part of a large dataset which  has already been the subject of many works in the literature. 
 See, for instance, Fokianos \textit{et al.} (2009), Fokianos and Neumann (2013), Davis and Liu (2016) (the series of July 2, 2002), 
 Doukhan and Kengne (2015) (the series of July 16, 2002), Diop and Kengne (2017) (the series of July 5, 2002) 
 and Br\"ann\"as and Quoreshi (2010). 
The data (the transaction during July 21, 2002) and its autocorrelation function displayed on Figure \ref{Graphe} (see (a) and (b)) show  three stylized facts: (i) a positive temporal dependence; (ii) the data are overdispersed (the empirical mean is $7.28$ while the empirical variance is $28.05$); (iii) presence of outliers is suspected.
Our purpose is to fit these data by taking into account a possible relationship between the number of transactions  and the volume-volatility. 
 This question has been investigated in several financial studies during the past two decades;  see, for example, Takaishi and Chen (2016), Belhaj \textit{et al.} (2015) and Louhichi (2011). 
In these works, the volume-volatility is found to exhibit a statistically significant impact on the trading volume (number of transactions or trade size).\\
Now, we consider the Poisson-INGARCH-X model given by
\begin{equation}\label{P_INGARCH_X(1,1)}
            	Y_{t}|\mathcal{F}_{t-1} \sim \mathcal{P}(\lambda_{t})~~\text{with}~~
            	 \lambda_{t}=  	\alpha_{0} + \alpha_{1} Y_{t-1} +  \beta \lambda_{t-i}                                                       +\gamma |V_{t-1}|, 
    \end{equation}  
    where
   $ \alpha_{0}>0$, $\alpha_{1}, \beta, \gamma  \geq 0$ and $V_t$ represents the volume-volatility at time $t$.
 Observe that if $\gamma=0$, then the model (\ref{P_INGARCH_X(1,1)}) reduces to the classical Poisson-INGARCH$(1,1)$ model.
  We first examine the adequacy of the fitted (based on the MLE) model with exogenous covariate by comparing it with the Poisson-INGARCH$(1,1)$ model (without exogenous covariate). 
  As evaluation criteria, we consider the estimated counterparts of the Pearson residuals defined by $e_t=(Y_t-  \lambda_t)/ \sqrt{\lambda_{t}}$. Under the true model, the process $\{ e_t\}$ is close to a white noise sequence with constant variance (see for instance, Kedem and Fokianos (2002)). 
 The comparison of these two models is based on the MSE of the Pearson residuals which is defined by $\sum_{t=1}^{n} e^2_t/(n-d)$,
  where $d$ denotes the number of the estimated parameters.  
 The approximated MSE of the Pearson residuals is $2.269$ for the model (\ref{P_INGARCH_X(1,1)}) and $2.298$ for the Poisson-INGARCH$(1,1)$;  which indicates a preference for the model with covariate.
 %
  %
The test of  the significance of the exogenous covariate  of Pedersen and Rahbek (2018), applied to the series also confirms these results. 
Under the model (\ref{P_INGARCH_X(1,1)}) (i.e, with $\gamma \neq 0$), the cumulative periodogram plot of the Pearson residuals is displayed in Figure \ref{Graphe}(d).
From this figure, the associated residuals appear to be uncorrelated over time.
This lends a substantial support to the choice of the model with exogenous covariate for  fitting these data. 

\medskip
 \noindent  
Since outliers are suspected in the data (see Figure \ref{Graphe}(a)), we apply the MDEPE to estimate the parameters of the model.    
To choose the optimal tuning parameter $\alpha$, we adopt the idea of Warwick and Jones (2005). It is based on the minimization of an Asymptotic approximation of the summed Mean Squared Error (AMSE) defined by 
\[
\widehat{\text{AMSE}} = (\widehat \theta_{\alpha,n}-\widehat \theta_{1,n})'(\widehat \theta_{\alpha,n}-\widehat \theta_{1,n})
+
\frac{1}{n} \text{Trace} \big[\widehat \Sigma_{\alpha,n}\big] ,
\]
 where $\widehat \theta_{1,n}$ is the MDPDE obtained with $\alpha=1$ and 
 $\widehat \Sigma_{\alpha,n}$ is a consistent estimator of the covariance matrix $\Sigma_\alpha$ given in Theorem \ref{th2}.
 To compute $\widehat \theta_{\alpha,n}$, the initial value $\widehat \lambda_{1}$ is set to be the empirical mean of the data and $\partial \widehat \lambda_{1}/\partial \theta $ is set to be the null vector.
 For different values of $\alpha$, the corresponding $\widehat{\text{AMSE}}$ are displayed in Table \ref{Res.appli}.
Based on the findings of this table, the optimal tuning parameter chosen is $\alpha=0.2$, which provides the minimum value of the $\widehat{\text{AMSE}}$ (indicated by the symbol $^{\bullet}$). 
Thus, the MDPDE  is more accurate than the MLE for this data.
%
%
With $\alpha=0.2$, the MDPDE applied to the  model (\ref{P_INGARCH_X(1,1)}) yields:
\begin{equation}\label{Estim_Appl}
     \widehat \lambda_{t}= \underset{(0.081)}{0.103}+ \underset{(0.024)}{0.144}Y_{t-1} 
     + \underset{(0.027)}{0.833}\widehat \lambda_{t-1}
     + \underset{(0.023)}{0.030}|V_{t-1}|,
\end{equation}
where in parentheses are the standard errors of the estimators obtained from the robust sandwich matrix. 
Figure \ref{Graphe}(c) displays the number of transactions ($Y_t$), the fitted values ($\widehat Y_t \coloneqq  \widehat \lambda_{t}$) and the $95\%$-prediction interval based on the underlying Poisson distribution. This figure shows that the fitted values capture reasonably well the dynamics of the observed process.
 \begin{table}[h!]
\scriptsize
\footnotesize
\centering
\caption{\it The $\widehat{\text{AMSE}} \times 10^2$ corresponding to some values of $\alpha$ in the model (\ref{P_INGARCH_X(1,1)}).}
\label{Res.appli}
\vspace{.2cm}
\begin{tabular}{ccccccccccccccc}

\Xhline{.9pt}
  &&&&&&&&&& \\

  $\alpha$&
             $0$&$0.05$&$0.1$&$0.15$&$0.2$&$0.25$&$0.3$&$0.35$&$0.4$&$0.45$&$0.5$&$0.75$&$1$\\
  \hline
  \rule[0cm]{0cm}{.4cm} 
  $\widehat{\text{AMSE}}$&
              $4.181$&$2.624$&$2.264$&$2.153$&$2.137^{\bullet}$&$2.156$&$2.205$&$2.282$&$2.374$&$2.488$&$2.624$&$3.628$&$5.267$\\    
  
       \Xhline{.9pt}
\end{tabular}
\end{table}

\section{Proofs of the main results}
 Without loss of generality, we only provide the proofs of Theorems \ref{th1} and \ref{th2} for $\alpha >0$.
 The proofs in the case of the maximum likelihood estimators (i.e., $\alpha=0$) can be done conventionally by using the classical methods.\\
\noindent Throughout the sequel, $C$ denotes a positive constant whom value may differ from an inequality to another.
 
 \subsection{Proof of Theorem \ref{th1}}
 
 We consider the following lemma.
 
 \begin{lem}\label{lem1} 
 Assume that the conditions of Theorem \ref{th1} hold. Then
 
\begin{equation}\label{eq_lem1}
\left\|\widehat H_{\alpha,n}(\theta) - H_{\alpha,n}(\theta)  \right\|_\Theta \limitepsn 0 .
 \end{equation}
 \end{lem}

\emph{\bf Proof of Lemma \ref{lem1}}\\
Remark that
\begin{align*}  
         \left\|\widehat H_{\alpha,n}(\theta) - H_{\alpha,n}(\theta)  \right\|_\Theta
         &\leq 
 \frac{1}{n} \sum_{t=1}^{n}\|\widehat{\ell}_{\alpha,t}(\theta)-\ell_{\alpha,t}(\theta) \|_\Theta\\    
 &\leq  I_{n,1}+I_{n,2},
 \end{align*} 
 where
 \begin{align*}
 I_{n,1}&=\frac{1}{n}\sum_{t=1}^{n}  \Big \|\sum_{y =0}^{\infty}\left\{ g (y|\widehat \eta_t(\theta))^{1+\alpha} -g (y|\eta_t(\theta))^{1+\alpha}  \right\}\Big\|_\Theta, \\
 I_{n,2}&=\big( 1+\frac{1}{\alpha}\big)\frac{1}{n} \sum_{t=1}^{n}\Big\| g (Y_t|\widehat \eta_t(\theta))^{\alpha} -g (Y_t|\eta_t(\theta))^{\alpha}  \Big\|_\Theta.
\end{align*} 

It suffices to show that $(i)$ $I_{n,1} \limitepsn 0$ and $(ii)$ $I_{n,2} \limitepsn 0$.  \\
\begin{enumerate}
	\item [$(i)$]
For any $t\geq 1$, we apply the mean value theorem at the function $\eta \mapsto \sum_{y =0}^{\infty}g (y| \eta)^{1+\alpha}$.
For any $\theta \in \Theta$, there exists $\tilde \eta_t(\theta)$ between $ \eta_t(\theta)$ and  $\widehat \eta_t(\theta)$ such that
\begin{align*}
\Big|\sum_{y =0}^{\infty}\left\{ g (y|\widehat \eta_t(\theta))^{1+\alpha} -g (y|\eta_t(\theta))^{1+\alpha}  \right\} \Big|
&\leq
(1+\alpha)\left| \widehat  \eta_t(\theta)-\eta_t(\theta)\right| \sum_{y =0}^{\infty} \Big| \frac{\partial g (y| \tilde \eta_t(\theta))}{\partial \eta}\Big|g (y|  \tilde \eta_t(\theta) )^{\alpha}\\
&\leq
C\left| \widehat  \eta_t(\theta)-\eta_t(\theta)\right| \sum_{y =0}^{\infty} \Big|\frac{\partial g (y| \tilde \eta_t(\theta))}{\partial \eta}\Big|\\ 
&\leq
C \big|  \eta(\widehat \lambda_t(\theta))-\eta(\lambda_t(\theta))\big|    \varphi ( \tilde \eta_t(\theta))\\
&\leq
C\big| \widehat \lambda_t(\theta)- \lambda_t(\theta) \big| \varphi ( \tilde \eta_t(\theta))  ~(\text{by virtue of } (\textbf{A6})).
%
\end{align*}
We deduce that 
 \begin{align*}
 I_{n,1}& \leq
 C\frac{1}{n}\sum_{t=1}^{n}  \big \| \big(\widehat  \lambda_t(\theta)-\lambda_t(\theta)\big) \varphi ( \tilde \eta_t(\theta)) \big\|_\Theta.
  \end{align*}  
  By using Kounias and Weng (1969), it suffices to show that
  \begin{equation}\label{proof_lem1_In_Kounias}
      \sum_{k \geq 1} \frac{1}{k}  \E \big \| \big(\widehat  \lambda_k(\theta)-\lambda_k(\theta)\big) \varphi ( \tilde \eta_k(\theta)) \big\|_\Theta< \infty.   
  \end{equation}  

%
%
From (\ref{App1_A0}) and the H\"{o}lder's inequality, we have
\begin{align*}
\E \big \| \big(\widehat  \lambda_k(\theta)-\lambda_k(\theta)\big) \varphi ( \tilde \eta_k(\theta)) \big\|_\Theta 
&\leq 
  \E \Big[     \big \| \varphi ( \tilde \eta_k(\theta)) \|_{\Theta}   \sum\limits_{\ell \geq k}  \alpha^{(0)}_{\ell} \left( Y_{k-\ell}+\left\|  X_{k-\ell}\right\|\right)\Big]\\
 &\leq \big(\E \big\| \varphi ( \tilde \eta_k(\theta)) \big\|^2_\Theta  \big)^{1/2} \Big(\E \Big[\big(\sum\limits_{\ell \geq k}  \alpha^{(0)}_{\ell,Y}  Y_{k-\ell}+ \sum\limits_{\ell \geq k}  \alpha^{(0)}_{\ell,X} \|  X_{k-\ell}\|\big)^2\Big]\Big)^{1/2}\\
 &\leq   \| \| \varphi ( \tilde \eta_k(\theta)) \|_\Theta \|_2  \Big[\sum\limits_{\ell \geq k}  \alpha^{(0)}_{\ell,Y}   \| Y_{k-\ell} \|_2 + \sum\limits_{\ell \geq k}  \alpha^{(0)}_{\ell,X}\| \| X_{k-\ell} \| \|_2 \Big]  \\
 & \leq   C \overline{\varphi}_k \sum\limits_{\ell \geq k} \big( \alpha^{(0)}_{\ell,Y} +\alpha^{(0)}_{\ell,X}\big)~ (\text{from (\textbf{A4}) and the stationary assumptions}) \\
 & \leq   C \frac{1}{k^{\gamma-1}} \overline{\varphi}_k ~(\text{from the Riemannian assumption }  (\ref{eq_th1})).  
\end{align*}
Hence,
\[ \sum_{k \geq 1} \frac{1}{k}  \E \big \| \big(\widehat  \lambda_k(\theta)-\lambda_k(\theta)\big) \varphi ( \tilde \eta_k(\theta)) \big\|_\Theta \leq C \sum_{k \geq 1} \frac{1}{k^{\gamma}} \overline{\varphi}_k < \infty   ~(\text{from }   (\ref{eq_th1})).\]
Therefore, (\ref{proof_lem1_In_Kounias}) holds and thus, $I_{n,1} \limitepsn 0$.

   \item [$(ii)$] By applying the mean value theorem at the function  $\eta \mapsto g (Y_t| \eta)^{\alpha}$ and from (\textbf{A6}), for all $\theta \in \Theta$, there exists $\tilde \eta_t(\theta)$ between $ \eta_t(\theta)$ and  $\widehat \eta_t(\theta)$ such that
      \begin{align*}
      \big| g (Y_t|\widehat \eta_t(\theta))^{\alpha} -g (Y_t|\eta_t(\theta))^{\alpha}   \big|
&\leq
(1+\alpha)\left| \widehat  \eta_t(\theta)-\eta_t(\theta)\right| \Big| \frac{\partial g (Y_t| \tilde \eta_t(\theta))}{\partial \eta}\Big|g (Y_t|  \tilde \eta_t(\theta) )^{\alpha-1}\\
&\leq
C\big| \widehat  \lambda_t(\theta)-\lambda_t(\theta)\big| \Big|\frac{1}{g (Y_t|  \tilde \eta_t )} \, \frac{\partial g (Y_t| \tilde \eta_t(\theta))}{\partial \eta} \Big|\\ 
&\leq
C\big| \widehat  \lambda_t(\theta)-\lambda_t(\theta)\big|\psi(\tilde \eta_t(\theta)).
     \end{align*}
 According to \cite{Kounias1969}, a sufficient condition for that $I_{n,2} \limitepsn 0$ is 
         \begin{equation}\label{proof_lem1_IIn_Kounias}
         \sum_{k \geq 1} \frac{1}{k} \E \Big [ \left\|\left(\widehat  \lambda_k(\theta)-\lambda_k(\theta)\right)\psi(\tilde \eta_k(\theta)) \right\|_\Theta \Big ] <\infty.
     \end{equation}
   According to (\ref{App1_A0}) and by using the same arguments as above, we get
  \begin{align*}
\E \big \| \big(\widehat  \lambda_k(\theta)-\lambda_k(\theta)\big) \psi( \tilde \eta_k(\theta)) \big\|_\Theta 
&\leq 
  \E \Big[     \big \| \psi ( \tilde \eta_k(\theta)) \|_{\Theta}   \sum\limits_{\ell \geq k}  \alpha^{(0)}_{\ell} \left( Y_{k-\ell}+\left\|  X_{k-\ell}\right\|\right)\Big]\\
 &\leq \big(\E \big\| \psi ( \tilde \eta_k(\theta)) \big\|^2_\Theta  \big)^{1/2} \Big(\E \Big[\big(\sum\limits_{\ell \geq k}  \alpha^{(0)}_{\ell} ( Y_{k-\ell}+\|  X_{k-\ell}\|)\big)^2\Big]\Big)^{1/2}\\
 &\leq   \| \| \psi( \tilde \eta_k(\theta)) \|_\Theta \|_2 \big[\sum\limits_{\ell \geq k}  \alpha^{(0)}_{\ell,Y}  \| Y_{k-\ell} \|_2 + \sum\limits_{\ell \geq k}  \alpha^{(0)}_{\ell,X} \| \| X_{k-\ell} \| \|_2  \big]  \\
 & \leq   C \overline{\psi}_k \sum\limits_{\ell \geq k} \big( \alpha^{(0)}_{\ell,Y}+\alpha^{(0)}_{\ell,X} \big)
  \leq   C \frac{1}{k^{\gamma-1}} \overline{\psi}_k .
\end{align*}
We deduce that
\[ \sum_{k \geq 1} \frac{1}{k}  \E \big \| \big(\widehat  \lambda_k(\theta)-\lambda_k(\theta)\big) \psi ( \tilde \eta_k(\theta)) \big\|_\Theta \leq C \sum_{k \geq 1} \frac{1}{k^{\gamma}} \overline{\psi}_k < \infty  \text{ from (\ref{eq_th1}) } .
\]
Hence, (\ref{proof_lem1_IIn_Kounias}) holds and thus, $I_{n,2} \limitepsn 0$. This achieves the proof of Lemma \ref{lem1}.
        $~~~~~~~~~~~~~~~~~~~~~~~~~~~~~~~~~~~~~~~~~~~~~~~~~~~~~~~~~~~~~ \blacksquare$\\
\end{enumerate}     
To complete the proof of Theorem \ref{th1}, we will show that: ({\bf 1.}) $\E \left[\left\| \ell_{\alpha,t}\right\|_\Theta\right] <\infty$ and ({\bf 2.}) the function $\theta \mapsto \E(\ell_{\alpha,1}(\theta))$ has a unique minimum at $\theta^*$.  
\begin{enumerate}
	\item [({\bf 1.})] For any $\theta \in \Theta$, we have
	\[ 
	\left| \ell_{\alpha,t}(\theta)\right| 
	 \leq \underset{y =0}{\overset{\infty}{\sum}} g (y|\eta_t(\theta))^{1+\alpha} +\big( 1+\frac{1}{\alpha}\big)g (Y_t|\eta_t(\theta))^{\alpha}
	 \leq \underset{y =0}{\overset{\infty}{\sum}} g (y|\eta_t(\theta)) +\big( 1+\frac{1}{\alpha}\big) 
	 \leq 2+\frac{1}{\alpha}.
	\]
	
	Hence, $\E \left[\left\| \ell_{\alpha,t}\right\|_\Theta\right] <\infty$.
	\item [({\bf 2.})] 
	Let $\theta \in \Theta$, with $\theta \neq \theta^*$.   We have
	\begin{align*}
	&\E(\ell_{\alpha,1}(\theta))-\E(\ell_{\alpha,1}(\theta^*))
	= \E \Big \{ \E \big[\big( \ell_{\alpha,1}(\theta)-\ell_{\alpha,1}(\theta^*)\big)\big|\mathcal{F}_{0} \big] \Big\}\\
	&= 
	\E \Bigg \{\sum_{y =0}^{\infty}  g (y|\eta_1(\theta))^{1+\alpha}
	-\sum_{y =0}^{\infty} g (y|\eta_1(\theta^*))^{1+\alpha}
	- \big( 1+\frac{1}{\alpha}\big)\E \bigg[g (Y_1|\eta_1(\theta))^{\alpha} -g (Y_1|\eta_1(\theta^*))^{\alpha}\Big|\mathcal{F}_{0}\bigg] \Bigg \}\\
	&= 
\E \Bigg \{\sum_{y =0}^{\infty}  g (y|\eta_1(\theta))^{1+\alpha}
	-\sum_{y =0}^{\infty} g (y|\eta_1(\theta^*))^{1+\alpha}
		- \big( 1+\frac{1}{\alpha}\big)\sum_{y =0}^{\infty} \bigg[\big(g (y|\eta_1(\theta))^{\alpha} -g (y|\eta_1(\theta^*))^{\alpha}\big)g (y|\eta_1(\theta^*))\bigg] \Bigg\}\\
		&= 
		\E \left\{ \sum_{y =0}^{\infty} \left[ g (y|\eta_1(\theta))^{1+\alpha}
	-g (y|\eta_1(\theta^*))^{1+\alpha}
		-\big( 1+\frac{1}{\alpha}\big)\big(g (y|\eta_1(\theta))^{\alpha}-g (y|\eta_1(\theta^*))^{\alpha}\big)g (y|\eta_1(\theta^*))\right] \right\}\\
		&= 
		\E \left\{ d_\alpha\big(g (\cdot|\eta_1(\theta)),g (\cdot|\eta_1(\theta^*))\big)\right\}\geq 0,
	\end{align*}
	where the equality holds a.s. if and only if 
	$\theta=\theta^*$ according to (\textbf{A0}), (\textbf{A3}), (\textbf{A6}).  Thus, the function $\theta \mapsto \E(\ell_{\alpha,1}(\theta))$ has a unique minimum at $\theta^*$. 
\end{enumerate}
Recall that since $\{(Y_{t},X_t),\,t\in \Z\}$ is stationary and ergodic, the process $\{\ell_{\alpha,t}(\theta),\,t\in \Z\}$ is also a stationary and ergodic sequence. Then, according to ({\bf 1.}), by the uniform strong law of large number applied on the process $\{\ell_{\alpha,t}(\theta),\,t\in \Z\}$, it holds that 
\begin{align}\label{eq1_proof_th1}
\left\| H_{\alpha,n}(\theta) - \E(\ell_{\alpha,1}(\theta)) \right\|_\Theta= \bigg\| \frac{1}{n} \sum_{t=1}^{n}\ell_{\alpha,t}(\theta) - \E(\ell_{\alpha,1}(\theta)) \bigg\|_\Theta \limitepsn 0 .
\end{align}
Hence, from Lemma \ref{lem1} and (\ref{eq1_proof_th1}), we have
\begin{align} \label{eq2_proof_th1}
\left\| \widehat H_{\alpha,n}(\theta) - \E(\ell_{\alpha,1}(\theta)) \right\|_\Theta 
&\leq 
\big\| \widehat H_{\alpha,n}(\theta) - H_{\alpha,n}(\theta)  \big\|_\Theta + 
\big\|  H_{\alpha,n}(\theta) - \E(\ell_{\alpha,1}(\theta)) \big\|_\Theta \limitepsn 0.
\end{align}
({\bf 2.}), (\ref{eq2_proof_th1}) and standard arguments lead to the conclusion. 
$~~~~~~~~~~~~~~~~~~~~~~~~~~~~~~~~~~~~~~~~~~~~~~~~~~~~~~~~~~~~~~~~~~~~~~~~~~~~~~~~~~~~~~~~ \blacksquare$\\

\subsection{Proof of Theorem \ref{th2}}
The following lemma is needed.

\begin{lem}\label{lem2} 
 Assume that the conditions of Theorem \ref{th2} hold. Then
 
\begin{equation}\label{eq_lem2}
\E\Big( \frac{1}{\sqrt{n}} \sum_{t=1}^{n}\Big\|\frac{ \partial \widehat \ell_{\alpha,t}(\theta)}{\partial \theta} - \frac{ \partial  \ell_{\alpha,t}(\theta)}{\partial \theta}  \Big\|_\Theta \Big) \limiten 0 .
 \end{equation}
 \end{lem}

\emph{\bf Proof of Lemma \ref{lem2}}\\
Remark that for all  $t \in \Z$,
\begin{align*}
 \frac{\partial \ell_{\alpha,t}(\theta)}{\partial \theta}
 &=
 (1+\alpha)\big[\underset{y =0}{\overset{\infty}{\sum}} \frac{\partial g (y|\eta_t(\theta))}{ \partial \theta}g (y|\eta_t(\theta))^{\alpha}  -\frac{\partial g (Y_t|\eta_t(\theta))}{ \partial \theta}g (Y_t|\eta_t(\theta))^{\alpha-1} \big] \\
& = (1+\alpha)\frac{\partial \eta_t(\theta)}{\partial \theta}\bigg[\underset{y =0}{\overset{\infty}{\sum}} \frac{\partial g (y|\eta_t(\theta))}{ \partial \eta}g (y|\eta_t(\theta))^{\alpha}  -\frac{\partial g (Y_t|\eta_t(\theta))}{ \partial \eta}g (Y_t|\eta_t(\theta))^{\alpha-1} \bigg] 
\\
& = (1+\alpha)\frac{\partial \lambda_t(\theta)}{\partial \theta}\eta'(\lambda_t(\theta))\bigg[\underset{y =0}{\overset{\infty}{\sum}} \frac{\partial g (y|\eta_t(\theta))}{ \partial \eta}g (y|\eta_t(\theta))^{\alpha}  -\frac{\partial g (Y_t|\eta_t(\theta))}{ \partial \eta}g (Y_t|\eta_t(\theta))^{\alpha-1} \bigg]
\\
&=h_\alpha(\lambda_t(\theta))\frac{\partial \lambda_t(\theta)}{\partial \theta},
\end{align*}
where the function $h_\alpha$ is defined in (\ref{def_h_alpha}). 
 $ \frac{\partial \widehat \ell_{\alpha,t}(\theta)}{\partial \theta}$ can be computed in the same way and by using the relation $|a_1b_1 - a_2b_2| \leq |a_1-a_2| |b_2| + |b_1-b_2| |a_1|, ~ \forall a_1, a_2, b_1, b_2 \in \R$,  we get,
\begin{align}\label{eq2_proof_lem2}
\Big\| \frac{\partial \widehat \ell_{\alpha,t}(\theta)}{\partial \theta} -\frac{\partial  \ell_{\alpha,t}(\theta)}{\partial \theta} \Big\|_\Theta
&= \Big\| h_\alpha(\widehat \lambda_t(\theta))\frac{\partial \widehat \lambda_t(\theta)}{\partial \theta} - h_\alpha( \lambda_t(\theta))\frac{\partial  \lambda_t(\theta)}{\partial \theta} \Big\|_\Theta
\nonumber\\
&\leq 
\big\|h_\alpha(\lambda_t(\theta))\big\|_\Theta \Big\| \frac{\partial \widehat \lambda_t(\theta)}{\partial \theta}-\frac{\partial  \lambda_t(\theta)}{\partial \theta} \Big\|_\Theta + \big\|h_\alpha(\widehat \lambda_t(\theta))-h_\alpha(\lambda_t(\theta))\big\|_\Theta  \Big\| \frac{\partial \widehat  \lambda_t(\theta)}{\partial \theta} \Big\|_\Theta . 
\end{align}
The mean value theorem applied to the function $ \lambda \mapsto h_\alpha(\lambda)$ gives,
\begin{align*}
\big|h_\alpha(\widehat \lambda_t(\theta))-h_\alpha(\lambda_t(\theta))\big| 
&=
\Big| \frac{\partial h_\alpha( \tilde \lambda_t( \theta))}{\partial \lambda}\Big| \big|\widehat \lambda_t(\theta)- \lambda_t(\theta) \big| \nonumber\\
&= \big| m_\alpha(\tilde \lambda_t(\theta))\big|\big|\widehat \lambda_t(\theta)- \lambda_t(\theta) \big|,
\end{align*}
where 
$ \tilde\lambda_t(\theta)$ is between $\widehat \lambda_t(\theta)$ and $\lambda_t(\theta)$; and the function $m_\alpha$ is defined in (\ref{def_m_alpha}).  
Thus,
\begin{multline}\label{lem2_E_partial_ell_bound_sum}
    \E\Big( \frac{1}{\sqrt{n}} \sum_{t=1}^{n}\Big\|\frac{ \partial \widehat \ell_{\alpha,t}(\theta)}{\partial \theta} - \frac{ \partial  \ell_{\alpha,t}(\theta)}{\partial \theta}  \Big\|_\Theta \Big)
 \leq
 \frac{1}{\sqrt{n}} \sum_{t=1}^{n} \E \Big[ \big\|h_\alpha(\lambda_t(\theta))\big\|_\Theta \Big\|\frac{\partial \widehat \lambda_t(\theta)}{\partial \theta}-\frac{\partial  \lambda_t(\theta)}{\partial \theta}\Big\|_\Theta \Big] \\
         + \frac{1}{\sqrt{n}} \sum_{t=1}^{n} \E \Big[ \big\|m_\alpha(\tilde \lambda_t(\theta))\big\|_\Theta \big\|\widehat \lambda_t(\theta)- \lambda_t(\theta) \big\|_\Theta \Big\| \frac{\partial \widehat  \lambda_t(\theta)}{\partial \theta} \Big\|_\Theta \Big].      
 \end{multline} 
From Assumption \textbf{A}$_1 (\Theta)$, for all $t \geq 1$, we have
    \begin{align}
    \Big \|\frac{\partial \widehat \lambda_t(\theta)}{\partial \theta}-\frac{\partial  \lambda_t(\theta)}{\partial \theta} \Big\|_\Theta
    \nonumber &=  \Big\| \frac{\partial}{\partial \theta}f_\theta ( Y_{t-1},\cdots,Y_{1},0,\cdots;X_{t-1},\cdots,X_{1},0,\cdots) - \frac{\partial}{\partial \theta}f_\theta  ( Y_{t-1},\cdots ;X_{t-1},\cdots ) \Big\|_\Theta \\
     \label{lem2_A1} &  \leq  \sum\limits_{\ell \geq t}  \alpha^{(1)}_{\ell,Y}  Y_{t-\ell}+ \sum\limits_{\ell \geq t}  \alpha^{(1)}_{\ell,X} \left\|  X_{t-\ell}\right\|,
      \end{align}
  and
 \begin{align}\label{lem2_A1_prime} 
    \Big \|\frac{\partial \widehat \lambda_t(\theta)}{\partial \theta} \Big\|_\Theta  & \leq  \Big\|\frac{\partial}{\partial \theta}f_\theta(0) \Big\|_\Theta  + \Big\| \frac{\partial}{\partial \theta}f_\theta ( Y_{t-1},\cdots,Y_{1},0,\cdots;X_{t-1},\cdots,X_{1},0,\cdots) - \frac{\partial}{\partial \theta}f_\theta(0) \Big\|_\Theta  \nonumber\\
     &  \leq \Big\|\frac{\partial}{\partial \theta}f_\theta(0) \Big\|_\Theta +\sum\limits_{\ell =1}^{t-1}   \alpha^{(1)}_{\ell,Y}  Y_{t-\ell}+ \sum\limits_{\ell =1}^{t-1}   \alpha^{(1)}_{\ell,X}\left\|  X_{t-\ell}\right\|  \nonumber \\
   & \leq C +\sum\limits_{\ell \geq1}   \alpha^{(1)}_{\ell,Y}  Y_{t-\ell}+\sum\limits_{\ell \geq1}   \alpha^{(1)}_{\ell,X}\left\|  X_{t-\ell}\right\|.     
      \end{align}     
  From (\ref{lem2_A1}), the H\"{o}lder's inequality and the assumption (\textbf{A7}), we get 
          \begin{align*}
    & \frac{1}{\sqrt{n}} \sum_{t=1}^{n} \E \Big[ \big\|h_\alpha(\lambda_t(\theta))\big\|_\Theta \Big\|\frac{\partial \widehat \lambda_t(\theta)}{\partial \theta}-\frac{\partial  \lambda_t(\theta)}{\partial \theta}\Big\|_\Theta \Big]\\  
    & \hspace{1cm}\leq \frac{1}{\sqrt{n}} \sum_{t=1}^{n}  \E \Big [ \big\|h_\alpha(\lambda_t(\theta))\big\|_\Theta \big(\sum\limits_{\ell \geq t}  \alpha^{(1)}_{\ell,Y}  Y_{t-\ell}+ \sum\limits_{\ell \geq t}  \alpha^{(1)}_{\ell,X} \left\|  X_{t-\ell}\right\| \big) \Big]\\
   &\hspace{1cm} \leq  \frac{1}{\sqrt{n}} \sum_{t=1}^{n} \| \|h_\alpha(\lambda_t(\theta)) \|_\Theta \|_2 
         \big[\sum\limits_{\ell \geq t}  \alpha^{(1)}_{\ell,Y} \|Y_{t-\ell}\|_2 + \sum\limits_{\ell \geq t}  \alpha^{(1)}_{\ell,X}\| \|  X_{t-\ell} \| \|_2\big] \\
 & \hspace{1cm} \leq  C \frac{1}{\sqrt{n}} \sum_{t=1}^{n} \sum\limits_{\ell \geq t}  \big(\alpha^{(1)}_{\ell,Y} +\alpha^{(1)}_{\ell,X} \big) 
   \leq  C \frac{1}{\sqrt{n}} \sum_{t=1}^{n} \frac{1}{t^{\gamma-1}} ~ (\text{from the condition }(\ref{eq_th2})) \\
 &\hspace{1cm} \leq  C \frac{1}{\sqrt{n}} (1+n^{2-\gamma}) \limiten 0, 
          \end{align*} 
 where the last convergence holds since $\gamma> 3/2$. 
 Therefore, the first term of the right hand side of (\ref{lem2_E_partial_ell_bound_sum}) converges to 0.
 Moreover, according to  (\ref{App1_A0}), (\ref{lem2_A1_prime}), the H\"{o}lder's inequality and the assumption (\textbf{A7}), we have     
   \begin{align*}
   & \frac{1}{\sqrt{n}} \sum_{t=1}^{n} \E \Big[ \|m_\alpha(\tilde \lambda_t(\theta)) \|_\Theta \big\|\widehat \lambda_t(\theta)- \lambda_t(\theta) \big\|_\Theta \Big\| \frac{\partial \widehat  \lambda_t(\theta)}{\partial \theta} \Big\|_\Theta \Big]   \\
  & \leq 
  \frac{1}{\sqrt{n}} \sum_{t=1}^{n} \E \Big[ \|m_\alpha(\tilde \lambda_t(\theta)) \|_\Theta \Big(  \sum\limits_{\ell \geq t}  \alpha^{(0)}_{\ell,Y} Y_{t-\ell}+\sum\limits_{\ell \geq t}  \alpha^{(0)}_{\ell,X}\left\|  X_{t-\ell}\right\|\Big)
   \Big( \Big\| \frac{\partial f_\theta(0)}{\partial \theta} \Big\|_\Theta +\sum\limits_{\ell =1}^{\infty} \alpha^{(1)}_{\ell,Y}  Y_{t-\ell}+\sum\limits_{\ell \geq 1}  \alpha^{(1)}_{\ell,X}\left\|  X_{t-\ell}\right\|  \Big) \Big ]\\
  & \leq 
  \frac{1}{\sqrt{n}} 
  \Big \| \Big \| m_\alpha(\tilde \lambda_t(\theta)) \Big\|_\Theta \Big\|_2 \cdot \Big\|\sum\limits_{\ell \geq t}  \alpha^{(0)}_{\ell,Y} Y_{t-\ell}+\sum\limits_{\ell \geq t}  \alpha^{(0)}_{\ell,X}\left\|  X_{t-\ell}\right\| \Big\|_4 \cdot 
         \Big\| \Big\| \frac{\partial f_\theta(0)}{\partial \theta} \Big\|_\Theta +\sum\limits_{\ell =1}^{\infty} \alpha^{(1)}_{\ell,Y}  Y_{t-\ell}+\sum\limits_{\ell \geq 1}  \alpha^{(1)}_{\ell,X}\left\|  X_{t-\ell}\right\| \Big\|_4  \\
 & \leq 
 \frac{1}{\sqrt{n}} \overline m_\alpha  \sum_{t=1}^{n} \Big( \sum\limits_{\ell \geq t}  \alpha^{(0)}_{\ell,Y} \|Y_{t-\ell} \|_4 +\sum\limits_{\ell \geq t}  \alpha^{(0)}_{\ell,X} \| \|X_{t-\ell} \| \|_4 \Big) 
  \Big( \Big\| \Big\| \frac{\partial f_\theta(0)}{\partial \theta} \Big\|_\Theta \Big\|_4 + \sum\limits_{\ell =1}^{\infty}  \alpha^{(1)}_{\ell,Y}   \|Y_{t-\ell} \|_4 + \sum\limits_{\ell =1}^{\infty}  \alpha^{(1)}_{\ell,X}\| \|X_{t-\ell} \| \|_4 \Big)\\
 & \leq  
 C \frac{1}{\sqrt{n}} \sum_{t=1}^{n}  \Big( \sum\limits_{\ell \geq t} \big(\alpha^{(0)}_{\ell,Y} +\alpha^{(0)}_{\ell,X} \big)\Big) \Big(C +C\big(\sum\limits_{\ell =1}^{\infty}  \alpha^{(1)}_{\ell,Y}+\sum\limits_{\ell =1}^{\infty}  \alpha^{(1)}_{\ell,X}\big)\Big)  \\
   & \leq 
    C \frac{1}{\sqrt{n}} \sum_{t=1}^{n}  \sum\limits_{\ell \geq t}  \big(\alpha^{(0)}_{\ell,Y} +\alpha^{(0)}_{\ell,X} \big)
   \leq  C \frac{1}{\sqrt{n}} \sum_{t=1}^{n} \dfrac{1}{t^{\gamma-1}} 
   \limiten 0 ~(\text{from the condition }(\ref{eq_th2}) \text{ and see above} ).
      \end{align*} 
Hence, the second term of the right hand side of (\ref{lem2_E_partial_ell_bound_sum}) converges to 0. This complete the proof of the lemma. $\blacksquare$ 
           ~~\\
            ~~\\

  The following lemma is also needed.
   \begin{lem}\label{lem3} 
 Assume that the conditions of Theorem \ref{th2} hold. Then
 
\begin{equation}\label{eq_lem3}
\bigg\|\frac{1}{n} \sum_{t=1}^{n}\frac{ \partial^2 \ell_{\alpha,t}(\theta)}{\partial \theta \partial \theta^T} - \E\Big( \frac{ \partial^2  \ell_{\alpha,1}(\theta)}{\partial \theta \partial \theta^T}\Big)  \bigg\|_\Theta \limitepsn 0 .
 \end{equation}
 \end{lem}

\emph{\bf Proof of Lemma \ref{lem3}}\\

Recall that $\frac{\partial \ell_{\alpha,t}(\theta)}{\partial \theta}=h_\alpha(\lambda_t(\theta))\frac{\partial \lambda_t(\theta)}{\partial \theta}$. Then, for $i,j \in \{1,\cdots, d\}$, we have 
\begin{align*}
\frac{ \partial^2 \ell_{\alpha,t}(\theta)}{\partial \theta_i \partial \theta_j}
&=\frac{\partial}{\partial \theta_j} \big(h_\alpha(\lambda_t(\theta))\frac{\partial \lambda_t(\theta)}{\partial \theta_i}\big)\\
&=
h_\alpha(\lambda_t(\theta))\frac{\partial^2 \lambda_t(\theta)}{\partial \theta_j \partial \theta_i}
+
\frac{\partial h_\alpha(\lambda_t(\theta))}{\partial \theta_j}\frac{\partial \lambda_t(\theta)}{\partial \theta_i} \\
&=
h_\alpha(\lambda_t(\theta))\frac{\partial^2 \lambda_t(\theta)}{\partial \theta_i \partial \theta_j}
+
\frac{\partial h_\alpha(\lambda_t(\theta))}{\partial \lambda}\frac{\partial \lambda_t(\theta)}{\partial \theta_i} \frac{\partial \lambda_t(\theta)}{\partial \theta_j}
\\
&=
h_\alpha(\lambda_t(\theta))\frac{\partial^2 \lambda_t(\theta)}{\partial \theta_i \partial \theta_j}
+
m_\alpha(\lambda_t(\theta))\frac{\partial \lambda_t(\theta)}{\partial \theta_i} \frac{\partial \lambda_t(\theta)}{\partial \theta_j}
.
\end{align*}
Let us show that $\E\Big( \left\|\frac{ \partial^2  \ell_{\alpha,t}(\theta)}{\partial \theta_i \partial \theta_j}\right\|_\Theta\Big) <\infty$, for all $i,j \in \{1,\cdots, d\}$.   
 From the  H\"{o}lder's inequality and the assumption (\textbf{A7}), we have
\begin{align}\label{eq1_proof_lem3}
\E\Big(\Big\|\frac{ \partial^2  \ell_{\alpha,t}(\theta)}{\partial \theta \partial_i \theta_j}\Big\|_\Theta\Big)\nonumber
&\leq
\left\|\left\|h_\alpha(\lambda_t(\theta))\right\|_\Theta\right\|_2 \Big\|\Big\|\frac{\partial^2 \lambda_t(\theta)}{\partial \theta_i \partial \theta_j}\Big\|_\Theta\Big\|_2
+
\left\|\|m_\alpha(\lambda_t(\theta))\right\|_\Theta\|_2 \cdot  \Big\|\Big\|\frac{\partial \lambda_t(\theta)}{\partial \theta_i}\Big\|_\Theta 
 \Big\|\frac{\partial \lambda_t(\theta)}{\partial \theta_j}\Big\|_\Theta\Big\|_2 \nonumber\\
 & \leq (\overline{h}_{\alpha}+1) \Big\|\Big\|\frac{\partial^2 \lambda_t(\theta)}{\partial \theta_i \partial \theta_j}\Big\|_\Theta\Big\|_2
+
\overline{m}_{\alpha} \Big\|\Big\|\frac{\partial \lambda_t(\theta)}{\partial \theta_i}\Big\|_\Theta\Big\|_4 \cdot 
 \Big\|\Big\|\frac{\partial \lambda_t(\theta)}{\partial \theta_j}\Big\|_\Theta\Big\|_4.
\end{align} 
Moreover, according to  the assumption \textbf{A}$_2 (\Theta)$, for all $t \in \Z$, we get
       
    \begin{align*}
    \Big \|\frac{\partial^2 \lambda_t(\theta)}{\partial \theta_i \partial \theta_j}\Big \|_\Theta
    \nonumber 
    & \leq 
    \Big\|\frac{\partial^2}{\partial \theta_i \partial \theta_j}f_\theta(0) \Big\|_\Theta
    +
     \Big\| \frac{\partial^2}{\partial \theta_i \partial \theta_j}f_\theta ( Y_{t-1},\cdots;X_{t-1},\cdots) - \frac{\partial^2}{\partial \theta_i \partial \theta_j}f_\theta(0) \Big\|_\Theta  \\
      \label{lem3_A2} &  \leq \Big\|\frac{\partial^2}{\partial \theta_i \partial \theta_j}f_\theta(0) \Big\|_\Theta +\sum\limits_{\ell =1}^{\infty} \alpha^{(2)}_{\ell,Y}  Y_{t-\ell}+ \sum\limits_{\ell =1}^{\infty} \alpha^{(2)}_{\ell,X}\left\|  X_{t-\ell}\right\|.     
      \end{align*} 
      Thus, according to  (\ref{eq1_proof_lem3}) and similar arguments  as in  (\ref{lem2_A1_prime}), we deduce  
      
      \begin{align*}
\E\Big(\Big\|\frac{ \partial^2  \ell_{\alpha,t}(\theta)}{\partial \theta \partial_i \theta_j}\Big\|_\Theta\Big)
&\leq
(\overline{h}_{\alpha}+1) \Big( \Big\|\Big\|\frac{\partial^2}{\partial \theta_i \partial \theta_j}f_\theta(0) \Big\|_\Theta\Big\|_2 +\sum\limits_{\ell =1}^{\infty}   \alpha^{(2)}_{\ell,Y}  \|Y_{t-\ell}\|_2+ \sum\limits_{\ell =1}^{\infty}   \alpha^{(2)}_{\ell,X} \left\|\left\|  X_{t-\ell}\right\|\right\|_2\Big)\\
& \hspace{3cm} +
\overline{m}_{\alpha} \Big( \Big\|\Big\|\frac{\partial}{\partial \theta}f_\theta(0) \Big\|_\Theta\Big\|_4 +\sum\limits_{\ell =1}^{\infty}   \alpha^{(1)}_{\ell,Y} \|Y_{t-\ell}\|_4+ \sum\limits_{\ell =1}^{\infty}\alpha^{(1)}_{\ell,X} \left\|\left\|  X_{t-\ell}\right\|\right\|_4\Big)^2\\
&\leq
(\overline{h}_{\alpha}+1)\Big( C +C\big(\sum\limits_{\ell =1}^{\infty}   \alpha^{(2)}_{\ell,Y} +\sum\limits_{\ell =1}^{\infty}   \alpha^{(2)}_{\ell,X} \big) \Big)
+
\overline{m}_{\alpha} \Big( C +C\big(\sum\limits_{\ell =1}^{\infty}   \alpha^{(1)}_{\ell,Y} +\sum\limits_{\ell =1}^{\infty}   \alpha^{(1)}_{\ell,X}\big) \Big)^2<\infty.
\end{align*}
%
Hence, from the  stationarity and ergodicity properties of the sequence $\big\{\frac{ \partial^2  \ell_{\alpha,t}(\theta)}{\partial \theta \partial \theta^T},\, t \in \Z \big\}$ and the uniform strong law of large numbers, it holds that 
\[
\Big\|\frac{1}{n} \sum_{t=1}^{n}\frac{ \partial^2 \ell_{\alpha,t}(\theta)}{\partial \theta \partial \theta^T} - \E\Big( \frac{ \partial^2  \ell_{\alpha,1}(\theta)}{\partial \theta \partial \theta^T}\Big)  \Big\|_\Theta \limitepsn 0 .
\]
Thus,  Lemma \ref{lem3} is verified.
           $~~~~~~~~~~~~~~~~~~~~~~~~~~~~~~~~~~~~~~~~~~~~~~~~~~~~~~~~~~~~~~~~~~~~~~~~~~~~~~~~~~~~~~~~~~~~ \blacksquare$\\
 ~\\
  
Now, we use the results of Lemma \ref{lem2} and \ref{lem3} to prove Theorem \ref{th1}.  
For $i \in \{ 1,\cdots,d\}$, by applying the Taylor expansion to the function $\theta \mapsto \frac{\partial}{\partial \theta_i} H_{\alpha,n}(\theta)$, there exists $\tilde \theta_{n,i}$ between $\widehat \theta_{\alpha,n}$ and $\theta^*$ such that
\[
\frac{\partial}{\partial \theta_i} H_{\alpha,n}(\widehat \theta_{\alpha,n})=\frac{\partial}{\partial \theta_i}H_{\alpha,n}(\theta^*)+ \frac{\partial^2}{\partial \theta \partial \theta_i} H_{\alpha,n}(\tilde \theta_{n,i}) (\widehat \theta_{\alpha,n}-\theta^*).
\]
    It comes that
%
\begin{align}\label{eq1_proof_th2}
    \sqrt{n}J_{n}(\widehat \theta_{\alpha,n})(\widehat \theta_{\alpha,n}-\theta^*)
    =
     \sqrt{n} \Big(\frac{\partial }{\partial \theta}H_{\alpha,n}(\theta^*)-\frac{\partial}{\partial \theta}H_{\alpha,n}(\widehat \theta_{\alpha,n})\Big),
    \end{align}    
    where
    \[
   J_{n}(\widehat \theta_{\alpha,n})= -\bigg(\frac{\partial^2}{\partial \theta \partial \theta_i} H_{\alpha,n}(\tilde \theta_{n,i}) \bigg)_{1\leq i \leq d}.
    \]
    We can rewrite (\ref{eq1_proof_th2}) as
    \[
     \sqrt{n}J_{n}(\widehat \theta_{\alpha,n})(\widehat \theta_{\alpha,n}-\theta^*)
    =
    \sqrt{n}\frac{\partial }{\partial \theta}H_{\alpha,n}(\theta^*)
     -\sqrt{n}\frac{\partial}{\partial \theta}\widehat H_{\alpha,n}(\widehat \theta_{\alpha,n})+ \sqrt{n}\Big(\frac{\partial}{\partial \theta}\widehat H_{\alpha,n}(\widehat \theta_{\alpha,n})-\frac{\partial}{\partial \theta}H_{\alpha,n}(\widehat \theta_{\alpha,n})\Big).
      \]
According to Lemma \ref{lem2}, it holds that
  \begin{align*}
  \E\Big(\sqrt{n}\Big|\frac{\partial}{\partial \theta}\widehat H_{\alpha,n}(\widehat \theta_{\alpha,n})-\frac{\partial}{\partial \theta}H_{\alpha,n}(\widehat \theta_{\alpha,n})\Big|\Big)
 & \leq 
  \E\Big(\sqrt{n} \Big\|\frac{\partial}{\partial \theta}\widehat H_{\alpha,n}(\theta)-\frac{\partial}{\partial \theta}H_{\alpha,n}(\theta)\Big\|_\Theta\Big)\\
 & \leq 
  \E\Big(\frac{1}{\sqrt{n}} \sum_{t=1}^{n}\bigg\|\frac{ \partial \widehat \ell_{\alpha,t}(\theta)}{\partial \theta} - \frac{ \partial  \ell_{\alpha,t}(\theta)}{\partial \theta}  \bigg\|_\Theta\Big) \limiten 0.
  \end{align*}
    Moreover, for $n$ large enough, $\frac{\partial}{\partial \theta}\widehat H_{\alpha,n}(\widehat \theta_{\alpha,n})=0$, since 
  $\widehat \theta_{\alpha,n}$ is a local minimum of the function $\theta \mapsto \widehat H_{\alpha,n}(\theta)$.\\
 So, for $n$ large enough, we have 
 \begin{align}\label{eq2_proof_th2}
    \sqrt{n}J_{n}(\widehat \theta_{\alpha,n})(\widehat \theta_{\alpha,n}-\theta^*)
    =
     \sqrt{n}\frac{\partial }{\partial \theta}H_{\alpha,n}(\theta^*) +o_P(1).
    \end{align}
    To complete the proof of Theorem \ref{th2}, we will show that
    
\begin{enumerate}
	\item [(i)] $\big\{\frac{ \partial  \ell_{\alpha,t}(\theta^*)}{\partial \theta}|\mathcal{F}_{t-1},\, t \in \Z \big\}$ is a stationary ergodic martingale difference sequence and $\E\big(\frac{ \partial  \ell_{\alpha,t}(\theta^*)}{\partial \theta}\big)^2<\infty$.
		\item [(ii)]  $J_{n}(\widehat \theta_{\alpha,n}) \limitepsn J_\alpha$.
			\item [(iii)] The matrix $J_\alpha$ is invertible.\\
\end{enumerate}

\begin{enumerate}
 \item [(i)] Recall that $\mathcal{F}_{t-1}=\sigma((Y_s,X_s),\, s \leq t-1)$ and $\frac{ \partial  \ell_{\alpha,t}(\theta^*)}{\partial \theta}=h_\alpha(\lambda_t(\theta^*))\frac{\partial \lambda_t(\theta^*)}{\partial \theta}$,
	where
	\[ 
	h_\alpha(\lambda_t(\theta^*))=(1+\alpha)\eta'(\lambda_t(\theta^*))\bigg[\underset{y =0}{\overset{\infty}{\sum}} \frac{\partial g (y|\eta_t(\theta^*))}{ \partial \eta}g (y|\eta_t(\theta^*))^{\alpha}  -\frac{\partial g (Y_t|\eta_t(\theta^*))}{ \partial \eta}g (Y_t|\eta_t(\theta^*))^{\alpha-1} \bigg].
	\]
	Since the functions $\lambda_t(\theta^*)$ and $\frac{\partial \lambda_t(\theta^*)}{\partial \theta}$ are $\mathcal{F}_{t-1}$-measurable, we have
	
 \[ \E\Big(\frac{ \partial  \ell_{\alpha,t}(\theta^*)}{\partial \theta}|\mathcal{F}_{t-1}\Big)  =\E\big(h_\alpha(\lambda_t(\theta^*))|\mathcal{F}_{t-1}\big)\frac{\partial \lambda_t(\theta^*)}{\partial \theta}  
	~  \text{ and } ~ 
	\E\big(h_\alpha(\lambda_t(\theta^*))|\mathcal{F}_{t-1}\big)= 0.
 \]
 Thus, $\E\Big(\frac{ \partial  \ell_{\alpha,t}(\theta^*)}{\partial \theta}|\mathcal{F}_{t-1}\Big) =0$. 
		Moreover, since $Y_t$, $\left\|  X_{t-\ell}\right\|$ and $\frac{\partial}{\partial \theta}f_\theta $ have $4th$-order moment, by using (\textbf{A7}) and the H\"{o}lder's inequality, we get 
		\begin{align*}
		\E\Big(\Big|\frac{ \partial  \ell_{\alpha,t}(\theta^*)}{\partial \theta}\Big|^2\Big)	&\leq 
		\E\Big(\left\| \left(h_\alpha(\lambda_t(\theta))\right)^2\right\|_\Theta \Big\|\frac{\partial \lambda_t(\theta)}{\partial \theta}\Big\|^2_\Theta\Big) \leq \left\|\left\| h_\alpha(\lambda_t(\theta))\right\|^2_\Theta\right\|_2  \Big\|\Big\|\frac{\partial \lambda_t(\theta)}{\partial \theta}\Big\|_\Theta\Big\|^2_4\\
	&\leq 
	\overline{h}_\alpha  \Big\|  \big\|\frac{\partial}{\partial \theta}f_\theta(0) \big\|_\Theta +\sum\limits_{\ell =1}^{\infty}   \alpha^{(1)}_{\ell,Y} Y_{t-\ell}+\sum\limits_{\ell =1}^{\infty}   \alpha^{(1)}_{\ell,X}\left\|  X_{t-\ell}\right\|\Big\|^2_4\\
	&\leq 
	\overline {h}_\alpha  \Big(  \big\|\big\|\frac{\partial}{\partial \theta}f_\theta(0) \big\|_\Theta\big\|_4 +\sum\limits_{\ell =1}^{\infty}   \alpha^{(1)}_{\ell,Y}  \left\|Y_{t-\ell}\right\|_4+\sum\limits_{\ell =1}^{\infty}   \alpha^{(1)}_{\ell,X}\left\|\left\|  X_{t-\ell}\right\|\right\|_4\Big)^2\\	
	&\leq 
	C \Big(  C +C\big(\sum\limits_{\ell =1}^{\infty}   \alpha^{(1)}_{\ell,Y}+\sum\limits_{\ell =1}^{\infty}   \alpha^{(1)}_{\ell,X} \big)\Big)^2	<\infty.
		\end{align*}
		 \item [(ii)] For any $j =1,\cdots,d$, we have
		 \begin{align*}
		& \Big| \frac{1}{n} \sum\limits_{t =1}^{n}\frac{\partial}{\partial \theta_j \partial \theta_i} \ell_{\alpha,t}(\tilde \theta_{n,i}) 
		 -\E\Big(\frac{\partial}{\partial \theta_j \partial \theta_i} \ell_{\alpha,1}( \theta^*)\Big)\Big| \\
		 &\leq
		 \Big| \frac{1}{n} \sum\limits_{t =1}^{n}\frac{\partial}{\partial \theta_j \partial \theta_i} \ell_{\alpha,t}(\tilde \theta_{n,i}) 
		 -\E\Big(\frac{\partial}{\partial \theta_j \partial \theta_i} \ell_{\alpha,1}(\tilde \theta_{n,i})\Big)\Big|
	 +
		 \Big|\E\Big(\frac{\partial}{\partial \theta_j \partial \theta_i} \ell_{\alpha,1}( \tilde \theta_{n,i})\Big)
		 -\E\Big(\frac{\partial}{\partial \theta_j \partial \theta_i} \ell_{\alpha,1}(\theta^*)\Big)\Big|\\
		 &\leq
		 \Big\| \frac{1}{n} \sum\limits_{t =1}^{n}\frac{\partial}{\partial \theta_j \partial \theta_i} \ell_{\alpha,t}(\theta) 
		 -\E\Big(\frac{\partial}{\partial \theta_j \partial \theta_i} \ell_{\alpha,1}( \theta)\Big) \Big\|_\Theta
		 +
		 \Big|\E\Big(\frac{\partial}{\partial \theta_j \partial \theta_i} \ell_{\alpha,1}( \tilde \theta_{n,i})\Big)
		 -\E\Big(\frac{\partial}{\partial \theta_j \partial \theta_i} \ell_{\alpha,1}(\theta^*)\Big)\Big|\\
		  &\limiten 0~~\text{(by virtue of Lemma \ref{lem3} and Theorem \ref{th1})}.
		 \end{align*}
This holds for any $i,j \in \{1,\cdots,d\}$. Thus,
\[
J_{n}(\widehat \theta_{\alpha,n})= -\bigg(\frac{\partial^2}{\partial \theta \partial \theta_i} H_{\alpha,n}(\tilde \theta_{n,i}) \bigg)_{1\leq i \leq d}
= -\bigg(\frac{1}{n} \sum\limits_{t =1}^{n}\frac{\partial}{\partial \theta \partial \theta_i} \ell_{\alpha,t}(\tilde \theta_{n,i})\bigg)_{1\leq i \leq d}
 \limitepsn -\E\Big(\frac{\partial}{\partial \theta \partial \theta^T} \ell_{\alpha,1}( \theta^*)\Big)=J_\alpha.
\]
\item [(iii)] Let $U$ be a non-zero vector of $\R^d$. We have 
\begin{align*}
U(-J_\alpha) U^T &=U\E\bigg( \E\Big(\frac{\partial}{\partial \theta \partial \theta^T} \ell_{\alpha,1}( \theta^*)|\mathcal{F}_{0}\Big)\bigg)U^T\\
&=\E\Big( \E\big(m_\alpha(\lambda_1(\theta^*))|\mathcal{F}_{0}\big)\Big(U \frac{\partial \lambda_t(\theta^*)}{\partial \theta}\Big) \Big(U\frac{\partial \lambda_t(\theta^*)}{\partial \theta}\Big)^T\Big), 
 \text{ since }  \E\big(h_\alpha(\lambda_1(\theta^*))|\mathcal{F}_{0}\big)= 0 \\
 &> 0
\end{align*}
according to
\begin{align*}
 \E\big(m_\alpha(\lambda_1(\theta^*))|\mathcal{F}_{0}\big)=
 (1+\alpha)\underset{y =0}{\overset{\infty}{\sum}} \bigg(\eta'(\lambda_1(\theta^*))\frac{\partial g (y|\eta(\lambda_1(\theta^*)))}{ \partial \eta}\bigg)^2 g (y|\eta(\lambda_1(\theta^*)))^{\alpha-1}>0,
\end{align*}
and the assumption (\textbf{A8}). 
This implies that the matrix $(-J_\alpha)$ is symmetric and positive definite. Thus, $J_\alpha$ is invertible. \\
\end{enumerate}
Now, from (i), we apply the central limit theorem for stationary ergodic martingale difference sequence.  
It follows that
\[
\sqrt{n}\frac{\partial }{\partial \theta}H_{\alpha,n}(\theta^*) =\frac{1}{\sqrt{n}}\sum\limits_{t =1}^{n}\frac{\partial }{\partial \theta}\ell_{\alpha,t}(\theta^*)
\limiteloin \mathcal{N}_d \left(0,I_\alpha \right),
\]
where
\[
I_\alpha= \E \Big[\Big(\frac{\partial }{\partial \theta}\ell_{\alpha,1}(\theta^*)\Big)\Big(\frac{\partial }{\partial \theta}\ell_{\alpha,1}(\theta^*)\Big)^T \Big].
\]
According to (\ref{eq2_proof_th2}), for $n$ large enough, (ii) and (iii) imply that 
\begin{align*}
    \sqrt{n}(\widehat \theta_{\alpha,n}-\theta^*)
    &=
     \left(J_{n}(\widehat \theta_{\alpha,n})\right)^{-1}\big(\sqrt{n}\frac{\partial }{\partial \theta}H_{\alpha,n}(\theta^*)\big) +o_P(1)\\
     &=
     J^{-1}_{\alpha}\big(\sqrt{n}\frac{\partial }{\partial \theta}H_{\alpha,n}(\theta^*)\big) +o_P(1)
     \limiteloin \mathcal{N}_d (0,J^{-1}_{\alpha} I_\alpha J^{-1}_{\alpha} ).
    \end{align*}
This completes the proof of Theorem \ref{th2}.\\
 $~~~~~~~~~~~~~~~~~~~~~~~~~~~~~~~~~~~~~~~~~~~~~~~~~~~~~~~~~~~~~~~~~~~~~~~~~~~~~~~~~~~~~~~~~~~~~~~~~~~~~~~~~~~~~~~~~~~~~~~~~~~~~~~~~~~~~~~~~~ \blacksquare$\\

\subsection{Proof of Proposition \ref{exp_existence}}
 Let $G_\lambda(y)$ be the cumulative distribution function of $g(y|\eta)$, with $\lambda = B(\eta)$, and its inverse $ G_\lambda^{-1}(u) := \inf\{ y \geq 0, ~ G_\lambda(y) \geq u \} $ for all $u \in [0,1]$.
 Let $(U_t)$ be a sequence of independent uniform $(0, 1)$ random variables and such that $\xi_t = (U_t,\varepsilon_t)$ is independent and identically distributed over time.
  Let us prove the existence of a $\tau$-weakly dependent stationary solution $(Y_t, \lambda_t, X_t)$ of (\ref{Model_Expo}) and (\ref{exp_Markov_cov}) satisfying  
  \begin{equation*} \label{proof_Y_t_G_U_t}
   Y_t=G_{\lambda_t}^{-1}(U_t).
  \end{equation*}
 %
%
Note that, by the Proposition A.2 of \cite{Davis2016}, we get for $\lambda_t$ and $\lambda'_t$,
\begin{equation} \label{proof_E_G_lambdat_diff}
 \E \big[ |G_{\lambda_t}^{-1}(U_t) - G_{\lambda'_t}^{-1}(U_t) | \big| \lambda_t, \lambda'_t  \big] = |\lambda_t - \lambda'_t|.
\end{equation}
For a solution $(Y_t, \lambda_t, X_t)$ that fulfills (\ref{Model_Expo}), (\ref{exp_Markov_cov}) and (\ref{proof_E_G_lambdat_diff}), set $Z_t = (Y_t, X_t)$.  We get
\begin{equation} \label{proof_Z_t_F}
 Z_t= (Y_t,  X_t) = \big(G_{\lambda_t}^{-1}(U_t),u(X_{t-1},\cdots ; \varepsilon_t) \big) := F(Z_{t-1},\cdots ; \xi_t) ~ \text{ with } ~ \lambda_t = f_\theta(Z_{t-1},\cdots).
\end{equation}
%
For a vector $z=(y,x) \in \N_0 \times \R^{d_x}$, define the norm $\| z \|_w = |y| + w_x \| x \| $ for some $ w_x > 0$.
According to Doukhan and Wintenberger (2008), it suffices to show that: (i)  $\E \| F(\pmb{z};\xi_t) \|_w < \infty $ for some $\pmb{z} \in \big(\N_0  \times \R^{d_x} \big)^\N$  and (ii) there exists a non-negative sequence $(\alpha_k(F))_{k \geq 1}$ satisfying $\sum_{k \geq 1} \alpha_k(F) < 1$ such that $\E \| F(\pmb{z};\xi_t) - F(\pmb{z}';\xi_t) \|_w \leq \sum_{k \geq 1} \alpha_k(F) \|z_k - z'_k \|_w $ for all $\pmb{z}, \pmb{z}' \in \big( \N_0 \times \R^{d_x} \big)^\N$.
 
 \medskip
 
 The point (i) holds directly from the assumptions on the functions $f_\theta$ and $u$.
 To prove (ii), for any $\pmb{z} \in \big( \N_0  \times \R^{d_x} \big)^\N$, we set $\pmb{z} = (\pmb{y}, \pmb{x} )$
with $\pmb{y} = (y_k)$, $\pmb{x} = (x_k)$ and  define $\lambda=f_\theta(z_1,\cdots)$.
For any $\pmb{z}, \pmb{z'} \in \big( \N_0 \times \R^{d_x} \big)^\N$, from (\ref{proof_E_G_lambdat_diff}), \textbf{A}$_0 (\Theta)$ and (\ref{exp_Lip_cov}), we have
\begin{align*}
 \E \| F(\pmb{z};\xi_t) - F(\pmb{z}';\xi_t) \|_w &= \E\big[ |G_{\lambda}^{-1}(U_t) - G_{\lambda'}^{-1}(U_t) | + w_x |u(\pmb{x}; \varepsilon_t) - u(\pmb{x}'; \varepsilon_t) |  \big] \\
 & \leq \sum_{k \geq 1} \alpha_{k,Y}^{(0)} |y_k - y_k'| + \sum_{k \geq 1} \alpha_{k,X}^{(0)}\| x_k - x_k' \|  + w_x \sum_{k \geq 1} \alpha_k(u) \|x_k - x'_k \| \\
 & \leq \sum_{k \geq 1} \alpha_{k,Y}^{(0)}|y_k - y_k'| + w_x \sum_{k \geq 1}   \big( \dfrac{\alpha_{k,X}^{(0)}}{w_x}  \| x_k - x_k' \| + \alpha_k(u) \|x_k - x'_k \| \big) \\
& \leq \sum_{k \geq 1} \alpha_k(F)\|z_k - z_k'\|_w 
\end{align*}
 with $\alpha_k(F) = \max \big\{ \alpha_{k,Y}^{(0)}, \frac{\alpha_{k,X}^{(0)}}{w_x} + \alpha_k(u) \big\}$. Thus, choose $w_x$  such that $ w_x > \frac{\sum_{k\geq 1} \alpha_{k,X}^{(0)}}{1- \sum_{k\geq 1} \max\{ \alpha_{k,Y}^{(0)}, \alpha_k(u)\}}$ to get $\sum_{k\geq 1} \alpha_k(F) < 1$.  
This complete the proof of the proposition.\\
$~~~~~~~~~~~~~~~~~~~~~~~~~~~~~~~~~~~~~~~~~~~~~~~~~~~~~~~~~~~~~~~~~~~~~~~~~~~~~~~~~~~~~~~~~~~~~~~~~~~~~~~~~~~~~~~~~~~~~~~~~~~~~~~~~~~~~~~~~~ \blacksquare$\\
 
%



\begin{thebibliography}{99}

\bibitem{Agosto2016}
 {\sc Agosto, A., Cavaliere, G., Kristensen, D. and Rahbek, A.}
 \newblock Modeling corporate defaults: Poisson autoregressions with exogenous covariates (PARX).
 \newblock {\em Journal of Empirical Finance 38(B)}, (2016),  640-663.

 \bibitem{Francq2016}
 {\sc Ahmad, A. and Francq, C.}
 \newblock Poisson QMLE of count time series models.
 \newblock {\em Journal of Time Series Analysis 37}, (2016),  291-314.
 
 \bibitem{Al-Osh1987}
 {\sc Al-Osh MA, Alzaid AA.}
 \newblock First order integer-valued autoregressive (INAR(1)) process.
 \newblock {\em Journal of Time Series Analysis 8}, (1987),  261-275.

 \bibitem{Al-Osh1990}
 {\sc Al-Osh MA, Alzaid AA.}
 \newblock Integer-valued $pth$-order autoregressive structure.
 \newblock {\em J. Appl. Probab. 27(2)}, (1990),  314-324.
 
 \bibitem{Aknouche2020}
 {\sc Aknouche A. and Francq  C.} 
 \newblock Count and duration time series with equal conditional
stochastic and mean orders. 
 \newblock {\em Econometric Theory}, (2020), 1-33.

  \bibitem{Basu1994}
{\sc Basu, S., Lindsay, B.G.}
 \newblock Minimum disparity estimation for continuous models: efficiency, distributions and robustness. 
\newblock {\em Ann. Inst. Statist. Math. 48}, (1994), 683-705.

 \bibitem{Basu1998}
 {\sc Basu, A., Harris, I.R., Hjort, N.L. and Jones, M.C.}
 \newblock Robust and efficient estimation by minimizing a density power divergence.
 \newblock {\em Biometrika 85}, (1998),  549-559.
 
 \bibitem{Belhaj2015}
 {\sc Belhaj, F., Abaoub, E. and Mahjoubi, M.}  
 \newblock Number of Transactions, Trade Size and the Volume-Volatility Relationship: An Interday and Intraday Analysis on the Tunisian Stock Market. 
 \newblock {\em International Business Research 8(6)}, (2015).
 
  \bibitem{Beran1977}
 {\sc Beran, R.}
  \newblock  Minimum Hellinger distance estimates for parametric models. 
  \newblock {\em Ann. Statist. 5}, (1977), 445-463.
  
  \bibitem{Brannas2010}
 {\sc Br\"ann\"as, K. and Shahiduzzaman Quoreshi, A.M.M.}
 \newblock Integer-valued moving average modelling of the number of transactions in stocks. 
 \newblock {\em Applied Financial Economics 20}, (2010), 1429-1440.
  
  \bibitem{Cox1981}
  {\sc Cox, D.R.}  
  \newblock Statistical analysis of time series: Some recent developments. 
   \newblock {\em Scandinavian Journal of Statistics 8}, (1981), 93-115.
 
  
  \bibitem{Cui2017}
  {\sc Cui, Y. and Zheng, Q.}  
  \newblock Conditional maximum likelihood estimation for a class of observation-driven time series models for count data. 
   \newblock {\em Statistics \& Probability Letters 123}, (2017), 193-201.
 

 
\bibitem{Davis2009}
{\sc Davis, R.A. and Wu, R.}
\newblock A negative binomial model for time series of counts.
\newblock {\em Biometrika 96}, (2009), 735-749.

 \bibitem{Davis2016}
 {\sc Davis, R.A. and Liu, H.}
 \newblock  Theory and Inference for a Class of Observation-Driven Models with Application to Time Series of Counts.
 \newblock {\em Statistica Sinica 26}, (2016), 1673-1707.
 
 \bibitem{Diop2017}
{\sc Diop, M.L. and Kengne, W.}
\newblock Testing parameter change in general integer-valued time series.
\newblock {\em J. Time Ser. Anal. 38}, (2017), 880-894.

\bibitem{Douc2017}
{\sc Douc, R., Fokianos, K., and Moulines, E.}
\newblock Asymptotic properties of quasi-maximum likelihood estimators in observation-driven time series models.
\newblock {\em Electronic Journal of Statistics  11}, (2017),  2707-2740.

\bibitem{Doukhan2012}
{\sc Doukhan, P., Fokianos, K., and Tjøstheim, D.}
\newblock On Weak Dependence Conditions for Poisson autoregressions.
\newblock {\em Statist. and Probab. Letters 82}, (2012), 942-948.

\bibitem{Doukhan2013}
{\sc Doukhan, P., Fokianos, K., and Tjøstheim, D.}
\newblock Correction to ''On weak dependence conditions for Poisson autoregressions'' [Statist. Probab. Lett. 82 (2012) 942-948].
\newblock {\em Statist. and Probab. Letters 83}, (2013), 1926-1927.


 \bibitem{Doukhan2015}
{\sc Doukhan, P. and Kengne, W. }
\newblock Inference and testing for structural change in general poisson autoregressive models.
\newblock {\em Electronic Journal of Statistics 9}, (2015), 1267-1314.
 
\bibitem{Doukhan2008}
{\sc Doukhan, P. and Wintenberger, O.}
\newblock Weakly dependent chains with infinite memory.
\newblock {\em Stochastic Process. Appl. 118},  (2008) 1997-2013.


 \bibitem{Ferland2006}
 {\sc Ferland, R., Latour, A. and Oraichi, D.}
 \newblock Integer-valued GARCH process.
 \newblock {\em J. Time Ser. Anal. 27}, (2006), 923-942.

 
 \bibitem{Fokianos2009}
 {\sc Fokianos, K., Rahbek, A. and Tjøstheim, D.}
 \newblock Poisson autoregression.
 \newblock {\em Journal of the American Statistical Association 104}, (2009), 1430-1439.
%

 \bibitem{Fokianos2010}
 {\sc Fokianos, K. and Fried, R.}
 \newblock Interventions in INGARCH processes.
 \newblock {\em J. Time Ser. Anal. 31}, (2010), 210-225.

 \bibitem{Fokianos2011}
{\sc Fokianos, K. and Tjøstheim, D.}
\newblock Log-linear Poisson autoregression.
\newblock {\em Journal of multivariate analysis 102},  (2011) 563-578.
%

 \bibitem{Fokianos2012_itv}
 {\sc Fokianos, K. and Fried, R.}
 \newblock Interventions in log-linear Poisson autoregression.
 \newblock {\em Statistical Modelling  12}, (2012), 1-24.


\bibitem{Fokianos2013b}
 {\sc Fokianos, K. and  Neumann, M.}
 \newblock A goodness-of-fit test for Poisson count processes.
 \newblock {\em Electronic Journal of Statistics 7}, (2013),  793-819.
 
 
\bibitem{Fokianos2019}
{\sc Fokianos, K. and Truquet, L. }
\newblock On categorical time series models with covariates.
\newblock {\em Stochastic processes and their applications 129}, (2019), 3446-3462.


\bibitem{Fokianos2020}
{\sc Fokianos, K., St{\o}ve, B., Tj{\o}stheim, D., and Doukhan, P. }
\newblock Multivariate count autoregression.
\newblock {\em Bernoulli 26}, (2020), 471-499.

\bibitem{Francq2019}
{\sc Francq, C. and  Thieu, L.Q.}  
\newblock QML inference for volatility models with covariates. 
{\em Econometric Theory 35}, (2019), 37-72.

 \bibitem{Fried2015}
{\sc Fried, R., Agueusop, I., Bornkamp, B., Fokianos, K., Fruth, J., Ickstadt, K.}
 \newblock Retrospective Bayesian outlier detection in INGARCH series. 
 \newblock {\em Statistics and Computing 25}, (2015), 365-374.



 
 

\bibitem{Kang2014}
{\sc Kang, J. and Lee, S.} 
\newblock Minimum density power divergence estimator for Poisson autoregressive models.
\newblock {\em Computational Statistics and Data Analysis 80}, (2014), 44-56.

\bibitem{Kedem2002}
 {\sc Kedem, B. and Fokianos, K.}
 \newblock Regression Models for Time Series Analysis.
 \newblock {\em Hoboken, Wiley, NJ}, (2002).
%

 \bibitem{Lee2017}
 {\sc  Kim, B. and Lee, S.}
 \newblock Robust estimation for zero-inflated Poisson autoregressive models based on density power divergence.
 \newblock {\em Journal of Statistical Computation and Simulation 87}, (2017), 2981-2996.
 

 \bibitem{Lee2019}
 {\sc  Kim B. and Lee, S.}
 \newblock Robust estimation for general integer-valued time series models.
 \newblock {\em Annals of the Institute of Statistical Mathematics}, (2019), 1-26.

 \bibitem{Kounias1969}
 {\sc Kounias, E.G. and Weng, T.-S }
 \newblock  An inequality and almost sure convergence.
 \newblock {\em Annals of Mathematical Statistics 33}, (1969), 1091-1093.


 \bibitem{Liboschik2017}
 {\sc  Liboschik, T., Fokianos, K. and Fried, R.}
 \newblock  tscount: An R package for analysis of count time series
following generalized linear models.
 \newblock {\em  Journal of Statistical Software 82}, (2017), 1-51.
 

\bibitem{Louhichi2011}
 {\sc Louhichi, W.}
 \newblock  What drives the volume-volatility relationship on Euronext Paris?
 \newblock {\em International Review of Financial Analysis 20}, (2011) 200-206.


 \bibitem{McKenzie1985}
 {\sc McKenzie, E.}
 \newblock  Some simple models for discrete variate time series.
 \newblock {\em Water Resour Bull. 21}, (1985), 645-650.
 
 %


 \bibitem{Pedersen2018}
 {\sc Pedersen, R.S. and Rahbek, A.}
 \newblock Testing Garch-X Type Models.
 \newblock {\em  Econometric Theory 35(5)}, (2018), 1-36.
 
 
\bibitem{Simpson1987}
 {\sc Simpson, D.G.}
  \newblock Minimum Hellinger distance estimation for the analysis of count data.
   \newblock {\em J. Amer. Statist. Assoc. 82}, (1987), 802-807.
   
 
 \bibitem{Takaishi2016}
 {\sc Takaishi, T., and Chen, T.T.}  
 \newblock The relationship between trading volumes, number of transactions, and stock volatility in GARCH models. 
 \newblock {\em Journal of Physics: Conference Series 738 (1)}, (2016).
 
 \bibitem{Tamura1986}
 {\sc Tamura, R.N. and Boos, D.D.} 
 \newblock Minimum Hellinger distance estimation for multivariate location and covariance. 
 \newblock {\em J. Amer. Statist. Assoc. 89}, (1986), 223-239.
 
 \bibitem{Warwick2005}
 {\sc Warwick, J. and Jones, M.C.}
 \newblock Choosing a robustness tuning parameter.
 \newblock {\em J. Statist. Comput. Simul. 75}, (2005), 581-588.
 



\end{thebibliography}
 \end{document}